\documentclass[11pt]{article}
\usepackage{amssymb,amsmath,amsopn,amsfonts}
\usepackage[dvips]{color}
\usepackage{colordvi,multicol}
\usepackage{epsf}
\usepackage{hyperref}
\usepackage{bbm}
\usepackage{xcolor}
\usepackage[T1]{fontenc} 

\newfam\msbfam
\font\tenmsb=msbm10 \textfont\msbfam=\tenmsb \font\sevenmsb=msbm7
\scriptfont\msbfam=\sevenmsb \font\fivemsb=msbm5
\scriptscriptfont\msbfam=\fivemsb

\def\cal{\mathcal}

\newcommand{\normmm}[1]{{\left\vert\kern-0.25ex\left\vert\kern-0.25ex\left\vert #1 
    \right\vert\kern-0.25ex\right\vert\kern-0.25ex\right\vert}}

\def\proof{\vspace{1mm}\noindent{\it Proof}\quad}

 \newtheorem{lemma}{Lemma}[section]
 \newtheorem{proposition}[lemma]{Proposition}
 \newtheorem{theorem}[lemma]{Theorem}
 \newtheorem{corollary}[lemma]{Corollary}
 \newtheorem{remark}[lemma]{Remark}
 \newtheorem{definition}[lemma]{Definition}

\numberwithin{equation}{section}

\def\bc{\begin{center}}
\def\ec{\end{center}}
\def\no{\noindent}

\def\f#1#2{\frac{#1}{#2}}

  \makeatletter
\newcommand{\rmnum}[1]{\romannumeral #1}
\newcommand{\Rmnum}[1]{\expandafter\@slowromancap\romannumeral #1@}
\makeatother

\begin{document}

\bigbreak
\title {{\bf Conservative stochastic 2-dimensional Cahn-Hilliard equation
\thanks{Research supported in part  by NSFC (No.11671035). Financial support by the DFG through the CRC 1283 "Taming uncertainty and profiting from randomness and low regularity in analysis, stochastics and their applications" is acknowledged.}\\}}
\author{{\bf Michael R\"{o}ckner}$^{\mbox{c,d},}$, {\bf Huanyu Yang}$^{\mbox{a,c,d}}$,{\bf Rongchan Zhu}$^{\mbox{b,c}}$,
\date {}
\thanks{E-mail address:  roeckner@math.uni-bielefeld.de(M. R\"{o}ckner), hyang@math.uni-bielefeld.de(H. Y. Yang), zhurongchan@126.com(R. C. Zhu),}\\ \\
\small $^{\mbox{a}}$School of Mathematical Science, University of Chinese Academy of Sciences, Beijing 100049, China,\\
\small $^{\mbox{b}}$Department of Mathematics, Beijing Institute of Technology, Beijing 100081, China,\\
\small $^{\mbox{c}}$ Department of Mathematics, University of Bielefeld, D-33615 Bielefeld, Germany,\\
\small $^{\mbox{d}}$ Academy of Mathematics and Systems Science, Chinese Academy of Sciences, Beijing 100190, China,}

\maketitle

\noindent {\bf Abstract}

We consider the stochastic 2-dimensional Cahn-Hilliard equation which is driven by the derivative in space of a space-time white noise. We use two different approaches to study this equation. First we prove that there exists a unique solution $Y$ to the shifted equation (\ref{1.3}). Then $X:=Y+{Z}$ is the unique solution to the stochastic Cahn-Hilliard equation, where ${Z}$ is the corresponding O-U process. Moreover, we use the Dirichlet form approach in \cite{Albeverio:1991hk} to construct a probabilistically weak solution to the original equation (\ref{1.1}) below. By clarifying the precise relation between the two solutions, we also get the restricted Markov uniqueness of the generator and the uniqueness of the martingale solutions to the equation (\ref{1.1}). Furthermore, we also obtain exponential ergodicity of the solutions.

 \no{\footnotesize{\bf Keywords}:\hspace{2mm}  stochastic quantization problem, Dirichlet forms, space-time white noise, Wick power, non-linear stochastic PDE}
\section{Introduction}
In this paper we show the well-posedness for the conservative stochastic Cahn-Hilliard equation
\begin{equation}\label{1.1}
    \left\{
   \begin{aligned}
   dX_t &=-\frac{1}{2}A\left(A X-:X^3:\right)dt+B dW_t,  \\
   X(0)&=z \in V_0^{-1},\\
   \end{aligned}
   \right.
 \end{equation}
on $\mathbb{T}^2$ in the probabilistically strong sense where $A=\Delta$, $B=\textrm{div}$. $W_t$ is an $L_0^2(\mathbb{T}^2,\mathbb{R}^2)$-cylindrical Wiener process, which is defined in Section \ref{s3}. $:X^3:$ denotes the Wick power, which is introduced in Section \ref{s3} and  the space $V_0^{-1}$ is defined similarly as the Sobolev space of order $-1$ (see Section \ref{s2}.
 
The Cahn-Hilliard equation is given by
$$
\partial_tu=-\Delta^2u-\Delta f(u),
$$
which was introduced by Cahn and Hilliard \cite{Cahn:1958ds} to study the phase separation of binary alloys. Here $f$ is the derivative of a free energy and generally $f$ is chosen as $f(u)=u^3-u$. The stochastic Cahn-Hilliard equation was first studied in \cite{Petschek:1983ik}, where Petschek and Metiu performed some numerical experiments for the stochastic Cahn-Hilliard equation driven by space-time white noise. In \cite{Elezovic:1991gy}, Elezovic and Mikelic proved the existence and uniqueness of a strong solution to the stochastic Cahn-Hilliard equation driven by trace-class noise. Then Da Prato and Debussche \cite{DaPrato:1996kk} proved existence and uniqueness of solutions for space-time white noise and obtained the existence and uniqueness of an invariant measure for trace-class noise. Later there many more papers appeared in which the authors study the properties of the solutions to the stochastic Cahn-Hilliard equations driven by trace-class noise (e.g. \cite{Debussche:2011ia,Anonymous:sW-sSAiV}).

For the conservative-type equation (\ref{1.1}), the Gibbs measure $\nu$ is formally given by
$$
\nu(d\phi)=c\exp\left(-\int_{\mathbb{T}^2} \frac{1}{4}:\phi^4:dx\right)\mu(d\phi) ,
$$
where $\mu$ is the Gaussian free field $\mathcal{N}(0,(-\Delta)^{-1})$, $c$ is a normalization constant, and $:\phi^4:$ is the fourth order Wick power of $\phi$. Thus the Gibbs measure $\nu$ is the restriction of  $\Phi_2^4$-field to the mass-conserving subspace, i.e. $\{f:\int f(x)dx=0\}$. Equation (\ref{1.1}) is sometimes called time-dependent Ginzburg-Landau (TDGL) equation of conservative type or Model B, while the stochastic Allen-Cahn equation (=dynamical $\Phi_2^4$-equation) is called the TDGL equation of non-conservative type or Model A (see \cite{HOHENBERG:1977vo,Funaki:2016iz}). { Since equation (\ref{1.1}) satisfies a conservation law, i.e. $\partial_t\int X(t,x)dx\equiv0$ for any solution $X$ to (\ref{1.1})} the class of equilibrium states is richer than for the stochastic Allen-Cahn equation.\footnote{Funaki's report}

In \cite{Parisi:1981wta} Parisi and Wu proposed a program for Euclidean quantum field theory based on getting Gibbs states of classical statistical mechanics as limiting distributions of stochastic processes, especially as solutions to non-linear stochastic differential equations. Then one can use the stochastic differential equations to study  properties of the Gibbs states. This procedure is called stochastic field quantization (see \cite{Administrator:2017uy}). The equation (\ref{1.1}) can be also viewed as a stochastic quantization equation for the conservative $\Phi_2^4$-field $\nu$.

Over the years, on the stochastic quantization of the $\Phi_2^4$-field, the  literature has kept on growing (see e.g. \cite{Administrator:2017uy,Albeverio:1991hk,Debussche:2003du,Mourrat:2015uo,Rockner:2017cca,Rockner:2015uh}) . The authors in these papers consider the following non-conservative stochastic quantization equation (Model A):
\begin{equation}\label{1.2}
dX_t=(A X-:X^3:)dt+dW_t.
\end{equation}
First results are due to Jona-Lasinio and Mitter \cite{Administrator:2017uy}. Using the Girsanov theorem, they constructed solutions to a modified equation on $\mathbb{T}^2$: 
\begin{equation}\label{1.2.1}
dX_t=(-\triangle+1)^{-\varepsilon}(\triangle X-:X^3:+aX)+(-\triangle+1)^{-\frac{\varepsilon}{2}}dW_t
\end{equation}
for $\frac{9}{10}<\varepsilon<1$. They also proved the ergodicity for (\ref{1.2.1}). In \cite{Albeverio:1991hk} Albeverio and R\"ockner studied (\ref{1.2}) using Dirichlet forms and constructed probabilistically weak solutions to (\ref{1.2}) for all $\varepsilon\in[0,1)$. In \cite{Mikulevicius:1999kt}, Mikulevicius and Rozovskii constructed martingale solutions to (\ref{1.2}) but the uniqueness remained open. In \cite{Debussche:2003du} Da Prato and Debussche considered the associated shifted equation to (\ref{1.2}) on $\mathbb{T}^2$ and proved local existence and uniqueness  of solutions in the probabilistically strong sense via a fixed point argument and then showed the non-explosion for almost every initial point by using the invariant measure. Recently Mourrat and Weber \cite{Mourrat:2015uo} showed the global existence and uniqueness for the shifted equation both on $\mathbb{T}^2$ and $\mathbb{R}^2$ for every initial point. Combining the results from the weak approach and strong approach, R\"ockner, Zhu and Zhu \cite{Rockner:2015uh} proved the restricted Markov uniqueness for the generator of (\ref{1.2}) and the uniqueness of the martingale problem to (\ref{1.2}) in \cite{Mikulevicius:1999kt} on $\mathbb{T}^2$ and $\mathbb{R}^2$. Furthermore, the ergodicity of (\ref{1.2}) on $\mathbb{T}^2$ has been obtained in \cite{Hairer:2016wl,Rockner:2017cca,Tsatsoulis:2016wm}.

For the conservative case, Funaki \cite{Funaki:1989va} proved the existence and uniqueness of solutions  to equation (\ref{1.1}) on $\mathbb{R}$ and in \cite{Debussche:2007kj} Debussche and Zambotti studied equation (\ref{1.1}) on $[0,1]$ with reflection. But for the higher dimensional case, even though the linear operator $\Delta^2$ gives much more regularity, the noise and hence the solutions are still so singular that the non-linear terms in (\ref{1.1}) are not well-defined in the classical sense. This difficulty is similar as in equation (\ref{1.2}). 

To overcome this difficulty, we use two approaches to study (\ref{1.1}). First we follow the idea in \cite{Debussche:2003du}, \cite{Mourrat:2015uo} and \cite{Rockner:2015uh} to split the solution to $X=Y+{Z}$, where ${Z}(t)=\int_{0}^t e^{-\frac{(t-s)}{2}A^2}BdW_s$.
Similarly as in the $\Phi_2^4$ case, $Y$ has better regularity than the solution to (\ref{1.1}) and satisfies the following \emph{shifted equation}:
  \begin{equation}\label{1.3}
     \left\{
   \begin{aligned}
   \frac{dY}{dt} &=-\frac{1}{2}A^2Y+\frac{1}{2}A\sum_{k=0}^{3}C_3^kY^{3-k}:{Z}^k: \\
   Y(0)&=z \\
   \end{aligned}
   \right.
  \end{equation}
  where $Z(t)=\int_0^te^{-\frac{t-s}{2}A^2}BdW_s$.
  In this paper we obtain the existence and uniqueness of the solution to (\ref{1.3}). The fixed point arguments for local well-posedness in \cite{Debussche:2003du} and \cite{Mourrat:2015uo}only hold for initial values in $\mathcal{C}^{-\frac{4}{3}+}$. Due to the singularity of the noise and the lack of a maximum principle and a uniform $L^p$-estimate, we only have a uniform $H^{-1}$-estimate (see Theorem \ref{t6.5}), which is not strong enough to combine it with local well-posedness (see Remark \ref{r4.7.2}). Instead, our argument is based on a classical compactness argument. We obtain the existence of global solutions starting from the uniform $H^{-1}$-estimate directly. Moreover we consider the solutions in $H^{-1}$ and use the $L^4$-integrability to obtain uniqueness for (\ref{1.3}). 
  
In addition, we use the method in \cite{Albeverio:1991hk} to construct the Dirichlet form for (\ref{1.1}) (see Theorem \ref{t6.2}), which is given by
$$
\mathcal{E}(\varphi,\psi) = \frac{1}{2}\int \langle \nabla \varphi,\nabla \psi\rangle _{V_0^{-1}}d\nu, \varphi,\psi\in \mathcal{F}C_b^{\infty},
$$  
where $\mathcal{F}C_b^{\infty}$ is defined in Section \ref{s5}. We note that the tangent space is chosen as $V_0^{-1}$ and the gradient operator $\nabla$ is also defined in $H^{-1}$. This is different from the Dirichlet form for (\ref{1.2}), where the tangent space is chosen as $L^2$ and the gradient is the $L^2$-derivative. By the integration by parts formula for $\nu$ we also obtain the closability for the bilinear form ($\mathcal{E}$,$\mathcal{F}C_b^{\infty}$). The closure ($\mathcal{E}$,$D(\mathcal{E})$) is a quasi-regular Dirichlet form, which enables us to construct a probabilistically weak solution to (\ref{1.1}). 
Then by clarifying the relation between this solution and the solution to (\ref{1.3}), we prove that $X-Z$, where $X$ is the solution obtained by the Dirichlet form approach, also satisfies the shifted equation (\ref{1.3}). It follows that the $\Phi_2^4$-field is an invariant measure for $X$. Then we obtain the Markov uniqueness in the restricted sense for the generator of the Dirichlet form restricted to $\mathcal{F}C_b^{\infty}$ and the uniqueness of probabilistically weak solutions to (\ref{1.1}) having $\nu$ as an invariant measure. 

We prove exponential ergodicity by two approaches. One simple and short way based on the Dirichlet form approach is presented in Remark \ref{r6.9}. Using a uniform estimate, an invariant measure can also be constructed by the Krylov-Bogoliubov method. We follow an idea from \cite{Tsatsoulis:2016wm} to prove the strong Feller property of the semigroup of the solution to the equation (\ref{1.1}). Then we obtain exponential convergence to the unique invariant measure of the semigroup  for every starting point.

Finally we comment on the motivations to study the Dirichlet form of equation (\ref{1.1}). First Dirichlet form theory plays an important role in the scaling limit of particle system. It was conjectured in  \cite{Giacomin:1999vn} that the stochastic Cahn-Hilliard equation is the scaling limit of the Kawasaki dynamics of the Ising-Kac model, while it has been proved in \cite{Bertini:1993ix,Mourrat:2017ww} that the stochastic Allen-Cahn equation is the scaling limit of Glauber dynamics of the Ising-Kac model. Until now, for the scaling limit of the Kawasaki dynamics, even in the $1-d$ case where no renormalization method involved, there is still no complete result (see \cite{IbertiThesis:X_zkJjfb}). To use Dirichlet form theory to identify the scaling limit of particle systems, the Markov uniqueness of the Dirichlet form is required. We hope that in future work we can use the restricted Markov uniqueness of the Dirichlet form obtained in this paper to study the scaling limit of the Kawasaki dynamics.
Another motivation to study this Dirichlet form is to
study spectral properties. As we have shown in Remark \ref{r6.9}, the spectral gap for the stochastic Cahn-Hilliard equation is controlled by the spectral gap of the stochastic Allen-Cahn equation in finite volume. But the situation is different from infinite volume case. On the whole space $\mathbb{R}^d$, it is expected that for the polynomial potential given by $\Phi^4-m\Phi^2$, the stochastic Allen-Cahn equation still has a spectral gap so that exponential ergodicity holds, while the stochastic Cahn-Hilliard equation looses this property. In fact, the Dirichlet form $(\Lambda,D(\Lambda))$ for equation (\ref{1.1}) on $\mathbb{R}^2$ can be directly constructed as the closure of the following bilinear form
\begin{equation}\label{dffs}
\Lambda(\varphi,\psi)=\int\langle\nabla f,\nabla g\rangle_{\dot{H}^{-1}}d\nu,\forall\varphi,\psi\in\mathcal{F}C_b^\infty(\dot{H}^{-1-}),
\end{equation}
where $\dot{H}^s$ is the homogeneous Sobolev space of order $s$, and $\nabla$ is the gradient in $\dot{H}^{-1}$. Similarly to Section \ref{s5}, it is easy to check that $(\Lambda,D(\Lambda))$ is quasi-regular and obtain a probabilistically weak solution directly. Since $(-\Delta)^{-1}$ is not bounded in $\dot{H}^{-1}$, the argument in Remark \ref{r6.9} fails in the case of $\mathbb{R}^2$. We hope to use (\ref{dffs}) to study spectral properties and functional inequalities of the stochastic Cahn-Hilliard equation on $\mathbb{R}^2$ in our future work. Moreover, from the viewpoint of particle systems, it has been proved in \cite{Bertini:1993ix,Mourrat:2017ww} that the stochastic Allen-Cahn equation is the scaling limit of the Glauber dynamics of the Ising-Kac model, while the stochastic Cahn-Hilliard equation is expected to be the scaling limit of the Kawasaki dynamics of the Ising-Kac model.The spectral gap for these two kinds of particle system was studied in \cite{Lu:1993vs}. The authors considered the model on a bounded domain with size $L$ and proved that as $L\to \infty$, the spectral gap for the Glauber dynamics remains strictly positive while the spectral gap for Kawasaki dynamics decays with a rate $L^{-2}$ (see also recent work for the continuum Sine-Gordon model \cite{Bauerschmidt:2019wb}).

 This paper is organized as follows: In Section \ref{s2} we collect some results related to Besov spaces. In Section \ref{s3} we study the solution to the linear equation and define the Wick power. In Section \ref{s4} we obtain the global existence and uniqueness of solutions to the shifted equation (\ref{1.3}). In Section \ref{s5} we obtain existence of probabilistically weak solutions via the Dirichlet form approach. By clarifying the relation between the two solutions we obtain that the $\Phi_2^4$-field $\nu$ is an invariant measure of $X$, Markov uniqueness in the restricted sense for the generator of the Dirichlet form restricted to $\mathcal{F}C_b^{\infty}$ and uniqueness of the probabilistically weak solutions to (\ref{1.1}). Finally we prove the strong Feller property and exponential ergodicity of the Markov semigroup associated to the solution to (\ref{1.1}) in Section \ref{s6}.
 
\section{Preliminaries}\label{s2}
Below we recall the definition of Besov spaces. For a general introduction to the theory of Besov spaces we refer to \cite{Bahouri:2011wc, Triebel:1978wb, Triebel:2006vl}. First we introduce the following
notations. Throughout the paper, we use the notation $a\lesssim b$ if there exists a constant $c > 0$
such that $a\leq cb$, and we write $a\backsimeq b$ if $a\lesssim b$ and $b\lesssim a$.
 The space of real valued infinitely differentiable functions of compact support is denoted by $\mathcal{D}(\mathbb{R}^d)$ or $\mathcal{D}$. The space of Schwartz functions is denoted by $\mathcal{S}(\mathbb{R}^d)$. Its dual, the space of tempered distributions, is denoted by $\mathcal{S}'(\mathbb{R}^d)$. The Fourier transform and the inverse Fourier transform are denoted by $\mathcal{F}$ and $\mathcal{F}^{-1}$, respectively.

 Let $\chi,\theta\in \mathcal{D}$ be non-negative radial functions on $\mathbb{R}^d$, such that

i. the support of $\chi$ is contained in a ball and the support of $\theta$ is contained in an annulus;

ii. $\chi(z)+\sum_{j\geq0}\theta(2^{-j}z)=1$ for all $z\in \mathbb{R}^d$.

iii. $\textrm{supp}(\chi)\cap \textrm{supp}(\theta(2^{-j}\cdot))=\emptyset$ for $j\geq1$ and $\textrm{supp}\theta(2^{-i}\cdot)\cap \textrm{supp}\theta(2^{-j}\cdot)=\emptyset$ for $|i-j|>1$.

We call such a pair $(\chi,\theta)$ dyadic partition of unity, and for the existence of dyadic partitions of unity we refer to \cite[Proposition 2.10]{Bahouri:2011wc}. The Littlewood-Paley blocks are now defined as
$$\Delta_{-1}u=\mathcal{F}^{-1}(\chi\mathcal{F}u)\quad \Delta_{j}u=\mathcal{F}^{-1}(\theta(2^{-j}\cdot)\mathcal{F}u).$$
 \vskip.10in

{\textbf{Besov spaces}}

For $\alpha\in\mathbb{R}$, $p,q\in [1,\infty]$, $u\in\mathcal{D}$ we define
$$\|u\|_{B^\alpha_{p,q}}:=(\sum_{j\geq-1}(2^{j\alpha}\|\Delta_ju\|_{L^p})^q)^{1/q},$$
with the usual interpretation as $l^\infty$ norm in case $q=\infty$. The Besov space $B^\alpha_{p,q}$ consists of the completion of $\mathcal{D}$ with respect to this norm and the H\"{o}lder-Besov space $\mathcal{C}^\alpha$ is given by $\mathcal{C}^\alpha(\mathbb{R}^d)=B^\alpha_{\infty,\infty}(\mathbb{R}^d)$. For $p,q\in [1,\infty)$,
$$B^\alpha_{p,q}(\mathbb{R}^d)=\{u\in\mathcal{S}'(\mathbb{R}^d):\|u\|_{B^\alpha_{p,q}}<\infty\}.$$
$$\mathcal{C}^\alpha(\mathbb{R}^d)\varsubsetneq \{u\in\mathcal{S}'(\mathbb{R}^d):\|u\|_{\mathcal{C}^\alpha(\mathbb{R}^d)}<\infty\}.$$
We point out that everything above and everything that follows can be applied to distributions on the torus (see \cite{Sickel:1985gg}, \cite{Stein:1972jf}). More precisely, let $\mathcal{S}'(\mathbb{T}^d)$ be the space of distributions on $\mathbb{T}^d$.  Besov spaces on the torus with general indices $p,q\in[1,\infty]$ are defined as
the completion of $C^\infty(\mathbb{T}^d)$ with respect to the norm $$\|u\|_{B^\alpha_{p,q}(\mathbb{T}^d)}:=(\sum_{j\geq-1}(2^{j\alpha}\|\Delta_ju\|_{L^p(\mathbb{T}^d)})^q)^{1/q},$$
and the H\"{o}lder-Besov space $\mathcal{C}^\alpha$ is given by $\mathcal{C}^\alpha=B^\alpha_{\infty,\infty}(\mathbb{T}^d)$.  We write $\|\cdot\|_{\alpha}$ instead of $\|\cdot\|_{B^\alpha_{\infty,\infty}(\mathbb{T}^d)}$ in the following for simplicity.  For $p,q\in[1,\infty)$, we have
$$B^\alpha_{p,q}(\mathbb{T}^d)=\{u\in\mathcal{S}'(\mathbb{T}^d):\|u\|_{B^\alpha_{p,q}(\mathbb{T}^d)}<\infty\}.$$
\begin{equation}\label{2.1}
\mathcal{C}^\alpha\varsubsetneq \{u\in\mathcal{S}'(\mathbb{T}^d):\|u\|_{\alpha}<\infty\}.
\end{equation}

Here we choose Besov spaces as  completions of smooth functions, which ensures that the Besov spaces are separable which has a lot of advantages for our analysis below.
 \vskip.10in

\textbf{Wavelet analysis}

We will also use wavelet analysis to determine the regularity of a distribution in a Besov space.
Below we briefly summarize wavelet analysis and refer to  work of Meyer \cite{Meyer:1995vf},
Daubechies \cite{Daubechies:1992ff} and \cite{Triebel:2006vl} for more details. For every $r > 0$, there exists a compactly supported function
$\varphi\in C^r(\mathbb{R})$ such that:

1. $\langle \varphi(\cdot),\varphi(\cdot-k)\rangle=\delta_{k,0}$ for every $ k\in \mathbb{Z}$;

2. There exist $\tilde{a}_k,k\in\mathbb{Z}$, with only finitely many non-zero values, and such that
$\varphi(x)=\sum_{k\in\mathbb{Z}}\tilde{a}_k\varphi(2x-k)$ for every $x\in\mathbb{R}$;

3. For every polynomial $P$ of degree at most $r$ and for every $x\in\mathbb{R}$,
$\sum_{k\in\mathbb{Z}}\int P(y)\varphi(y-k)dy \varphi(x-k) = P(x)$.

Given such a function $\varphi$, we define for every $ x\in \mathbb{R}^d$ the recentered and
rescaled function $\varphi^n_x$ as follows
$$\varphi^n_x(y) := \Pi_{i=1}^d2^{\frac{n}{2}}\varphi(2^n(y_i-x_i)).$$
Observe that this rescaling preserves the $L^2$-norm. We let $V_n$ be the subspace of
$L^2(\mathbb{R}^d)$ generated by $\{\varphi^n_x: x\in \Lambda_n\}$, where
$$\Lambda_n := \{( 2^{-n}k_1,...,2^{-n}k_d) : k_i\in\mathbb{Z}\}.$$
An important property of wavelets is the existence of a finite set $\Psi$	 of compactly
supported functions in $C^r$ such that, for every $n\geq0$, the orthogonal complement
of $V_n$ inside $V_{n+1}$ is given by the linear span of all the $\psi^n_x, x\in \Lambda_n, \psi\in \Psi$.
 For every $n\geq0$
$$\{\varphi^n_x, x\in\Lambda_n\}\cup \{\psi^m_x: m\geq n, \psi\in\Psi, x\in \Lambda_m\},$$
forms an orthonormal basis of $L^2(\mathbb{R}^d)$.
This wavelet analysis allows one to identify a countable collection of conditions
that determine the regularity of a distribution.

 Setting ${\Psi}_\star={\Psi}\cup\{\varphi\}$, by results on weighted Besov spaces (see \cite[(2.2), (2.3), (2.4)]{Rockner:2015uh} and its reference for details),  we know that for $p\in(1,\infty)$, $\alpha\in\mathbb{R}$, $f\in {\mathcal{C}}^{\alpha}$
\begin{equation}\label{2.4}
\|f\|^p_{\alpha}\lesssim\sum_{n=0}^\infty 2^{n(\alpha+1)p}\sum_{\psi\in\Psi_\star}\sum_{x\in\Lambda_n}|\langle f,\psi^n_{x}\rangle|^pw(x)^p.
\end{equation}
where $w(x)=(1+|x|^2)^{-\frac{\sigma}{2}}, \sigma>0$ .

 \vskip.10in

\textbf{Estimates on the torus}

In this part we give estimates on the torus for later use.
 Set $\Lambda= (I-\Delta)^{\frac{1}{2}}$. For $s\geq0, p\in [1,+\infty]$ we use $H^{s}_p$ to denote the subspace of $L^p(\mathbb{T}^d)$, consisting of all  $f$   which can be written in the form $f=\Lambda^{-s}g, g\in L^p(\mathbb{T}^d)$ and the $H^{s}_p$ norm of $f$ is defined to be the $L^p$ norm of $g$, i.e. $\|f\|_{H^{s}_p}:=\|\Lambda^s f\|_{L^p(\mathbb{T}^d)}$.

\vskip.10in
To study (1.1) in the finite volume case, we need several important properties of Besov spaces on the torus and we recall the following Besov embedding theorems on the torus first (c.f. \cite[Theorem 4.6.1]{Triebel:1978wb}, \cite[Lemma A.2]{Gubinelli:2015jw}, \cite[Remark 3, Section 2.3.2]{Triebel:1992dx}):
\vskip.10in
 \begin{lemma}\label{l2.1}
  (\rmnum{1}) Let $1\leq p_1\leq p_2\leq\infty$ and $1\leq q_1\leq q_2\leq\infty$, and let $\alpha\in\mathbb{R}$. Then $B^\alpha_{p_1,q_1}(\mathbb{T}^d)$ is continuously embedded in $B^{\alpha-d(1/p_1-1/p_2)}_{p_2,q_2}(\mathbb{T}^d)$.

 (\rmnum{2}) Let $s\geq0$, $1<p<\infty$, $\varepsilon>0$. Then
 $ H^{s+\varepsilon}_p\subset B^{s}_{p,1}(\mathbb{T}^d)\subset B^{s}_{1,1}(\mathbb{T}^d)$.

 (\rmnum{3}) Let $1\leq p_1\leq p_2<\infty$ and let $\alpha\in\mathbb{R}$. Then $H^\alpha_{p_1}\subset H^{\alpha-d(1/p_1-1/p_2)}_{p_2}$.
 
 (\rmnum{4}) Let $0<q\leq\infty$, $1\leq p\leq\infty$ and $s>0$. Then $B_{p,q}^s\subset L^p$.
\end{lemma}
Here  $\subset$ means that the embedding is continuous and dense.

\vskip.10in

We recall the following Schauder estimates, i.e. the smoothing effect of the heat semigroup, for later use.

\vskip.10in
\begin{lemma}(\cite[Lemma A.7]{Gubinelli:2015jw}) \label{l2.2}
 Let $u\in B^{\alpha}_{p,q}(\mathbb{T}^d)$ for some $\alpha\in \mathbb{R}, p,q\in [1,\infty]$. Then for every $\delta\geq0$
$$\|e^{-tA^2}u\|_{B^{\alpha+\delta}_{p,q}(\mathbb{T}^d)}\lesssim t^{-\delta/4}\|u\|_{B^{\alpha}_{p,q}(\mathbb{T}^d)}.$$
\end{lemma}

\vskip.10in
One can  extend the multiplication on suitable Besov spaces and also  have the duality properties of Besov spaces from \cite[Chapter 4]{Triebel:1978wb}:
\vskip.10in

\begin{lemma} \label{l2.3}
(i) The bilinear map $(u; v)\mapsto uv$
extends to a unique continuous map from $\mathcal{C}^\alpha\times \mathcal{C}^\beta$ to $\mathcal{C}^{\alpha\wedge\beta}$ if and only if $\alpha+\beta>0$.

(ii) Let $\alpha\in (0,1)$, $p,q\in[1,\infty]$, $p'$ and $q'$ be their conjugate exponents, respectively. Then the mapping  $(u; v)\mapsto \int uvdx$  extends to a unique continuous bilinear form on $B^\alpha_{p,q}(\mathbb{T}^d)\times B^{-\alpha}_{p',q'}(\mathbb{T}^d)$.
\end{lemma}

\vskip.10in
We recall the following interpolation and  multiplicative inequalities for the elements in $H^s_p$, which is required for the a-priori estimate in Section \ref{s4} (cf. \cite[Theorem 4.3.1]{Triebel:1978wb}, \cite[Lemma 2.1]{Rockner:2015ej}, \cite[Theorem 2.80]{Bahouri:2011wc}):
 \vskip.10in

\begin{lemma}\label{l2.4}
 (\rmnum{1})  Suppose that $s\in (0,1)$ and $p\in (1,\infty)$. Then for $u\in H^1_p$
$$\|u\|_{H^s_p}\lesssim \|u\|_{L^p(\mathbb{T}^d)}^{1-s}\|u\|_{H^1_p}^s.$$

(\rmnum{2}) Suppose that $s>0$ and $p\in (1,\infty)$. If $u,v\in C^{\infty}(\mathbb{T}^2)$ then
$$\|\Lambda^s(uv)\|_{L^p(\mathbb{T}^d)}\lesssim\|u\|_{L^{p_1}(\mathbb{T}^d)}\|\Lambda^sv\|_{L^{p_2}(\mathbb{T}^d)}+\|v\|_{L^{p_3}(\mathbb{T}^d)}
\|\Lambda^su\|_{L^{p_4}(\mathbb{T}^d)},$$
with $p_i\in (1,\infty], i=1,...,4$ such that
$$\frac{1}{p}=\frac{1}{p_1}+\frac{1}{p_2}=\frac{1}{p_3}+\frac{1}{p_4}.$$

(\rmnum{3}) Suppose that $s_1<s_2$ and $1\leq p,q\leq\infty$. Then for $u\in B_{p,q}^{s_2}$ and $\forall \theta\in(0,1)$
$$
||u||_{B_{p,q}^{\theta s_1+(1-\theta)s_2}}\leq ||u||_{B_{p,q}^{s_1}}^{\theta}||u||_{B_{p,q}^{s_2}}^{1-\theta}.
$$
\end{lemma} 
\vskip.10in

We also collect some important properties for the multiplicative structure of Besov spaces from \cite{Mourrat:2015uo} and \cite{Triebel:2006vl}.
\vskip.10in
\begin{lemma}(\cite[Corollary 3.19, Corollary 3.21]{Mourrat:2015uo}) \label{l2.5}
(1) For $\alpha>0, p_1, p_2, p, q\in[1,\infty],  \frac{1}{p_1}+\frac{1}{p_2}=\frac{1}{p}$, the bilinear map $(u; v)\mapsto uv$
extends to a unique continuous bilinear map from $B^{\alpha}_{p_1,q}\times B^{\alpha}_{p_2,q}$ to $B^{\alpha,}_{p,q}$.

(2) For $\alpha<0, \alpha+\beta>0,  p_1, p_2, p, q\in[1,\infty],   \frac{1}{p_1}+\frac{1}{p_2}=\frac{1}{p}$, the bilinear map $(u; v)\mapsto uv$
extends to a unique continuous bilinear map from $B^{\alpha}_{p_1,q}\times B^{\beta}_{p_2,q}$ to $B^{\alpha}_{p,q}$.

\end{lemma}
\vskip.10in
 
\textbf{Notations} 

  Let $L$ denote the space $L^2(\mathbb{T}^2)=[0,1]^2$, where $\mathbb{T}^2$ is the 2 dimensional torus and we use $\langle \cdot,\cdot\rangle$ to denote the inner product in $L$. $A$ is the Laplacian operator on $L$, that is,
  \begin{equation}\label{3.1}
   D(A)=H^2_2(\mathbb{T}^2), A=\frac{\partial^2}{\partial x^2}+\frac{\partial^2}{\partial y^2}.
  \end{equation}
  $A$ is a self-adjoint operator in $L$, with complete orthonormal system $(e_n)_n$ of eigenvectors in $L$, given by
  \begin{align*}
 & e_0(x):=1,e_{(k_1,0)}(x)=\sqrt{2}e^{i\pi k_1x_1}, e_{(0,k_2)}(x)=\sqrt{2}e^{i\pi k_2x_2},\\ &e_k(x):=2e^{i\pi (k_1x_1+k_2x_2)},k_1k_2\neq 0.
  \end{align*}
  Then we have $Ae_k=-\lambda_k e_k$, where $\lambda_k=|k|^2{\pi}^2,k=(k_1,k_2)\in \mathbb{Z}^2,|k|^2=k_1^2+k_2^2$.  We also introduce a notation for the average of $h\in \mathcal{S}'(\mathbb{T}^2)$:
  $$m(h):={}_{\mathcal{S}'}\langle h,e_0\rangle_{\mathcal{S}}.$$
  
  For any $\alpha\in \mathbb{R}$, we define
  $$
V^{\alpha}:=\{u\in \mathcal{S}':\sum_{k}\lambda_k^{\alpha}|{}_{\mathcal{S}'}\langle u,e_k \rangle_{\mathcal{S}}|^2<\infty\}.  
  $$
For any $u,v\in V^{\alpha}$, define 
$$
\langle u,v \rangle_{V^{\alpha}}:=m(u){m(v)}+\sum_{k}\lambda_k^{\alpha}{}_{\mathcal{S}'}\langle u,e_k \rangle_{\mathcal{S}}{{}_{\mathcal{S}'}\langle v,e_k \rangle_{\mathcal{S}}}.
$$
It's easy to see that $(V^{\alpha},\langle \cdot,\cdot \rangle_{V^{\alpha}})$ is a Hilbert space and $V^{\alpha}\simeq H^{\alpha}_2$. Then for any $s,\alpha\in \mathbb{R}$, we can define a bounded linear operator $(-A)^{s}:V^\alpha\to V^{\alpha-2s}$ by:
$$(-A)^{s}u=\sum_{k\in\mathbb{Z}^2\setminus\{(0,0)\}}\lambda_k^{s}u_ke_k,$$
where $u=\sum_{k}u_ke_k\in V^{\alpha}$.
In particular, we set $Q:=(-A)^{-1}$ and extend it to a one-to-one bounded linear operator $\bar{Q}$ by
  \begin{equation}\label{3.2}
   \bar{Q}h:=Qh+m(h)e_0.
  \end{equation}
  Note that 
  \begin{equation}\label{3.3}
   Qe_k=
    \begin{cases}
     \frac{1}{\lambda_k}e_k & k\neq (0,0), \\
     0 & k=(0,0),
    \end{cases} 
  \end{equation}
  and
  \begin{equation}\label{3.4}
   \bar{Q}e_k=
    \begin{cases}
     \frac{1}{\lambda_k}e_k & k\neq (0,0) ,\\
     e_0 & k=(0,0).
    \end{cases} 
  \end{equation}
  Then we have
  $$\langle u,v\rangle _{V^{\alpha}}:=\langle \bar{Q}^{-\alpha/2}u,\bar{Q}^{-\alpha/2}v\rangle,$$
  and $\bar{Q}^{s}: V^{\alpha}\to V^{\alpha+2s}$ is an isomorphism for any $\alpha,s\in\mathbb{R}$, since
  $$
  \langle \bar{Q}^{s}u,\bar{Q}^{s}v\rangle _{V^{\alpha+2s}}=\langle u,v\rangle _{V^{\alpha}}.
  $$
  
 We also set
  $$V_0^{\alpha}:=\{h \in V^{\alpha}:\langle h,e_0\rangle _{V^{\alpha}}=0\},$$
  and denote $L_0^2:=V_0^0$. Let $\Pi$ denote the symmetric projector of $V^{\alpha}$ on $V_0^{\alpha}$, that is,
\begin{equation}\label{2.7.1}
\Pi:V^{\alpha} \to V_0^{\alpha},\Pi h:=h-m(h).
\end{equation}

Moreover, we define 
$$V^{\alpha}(\mathbb{T}^2,\mathbb{R}^2):=\{f=(f_1,f_2):f_i\in V^{\alpha},i=1,2\},$$
and similarly
$$
V_0^{\alpha}(\mathbb{T}^2,\mathbb{R}^2):=\{f=(f_1,f_2):f_i\in V_0^{\alpha},i=1,2\}.
$$
 
In this paper, we consider the initial value  and the reference measure on $V_0^{\alpha}$ for simplicity. For the general case, we refer to \cite{Debussche:2007kj}.
 
  \section{The Linear Equation and Wick Powers}\label{s3}
  
  We consider the O-U process
  \begin{equation}\label{4.1}
    \left\{
   \begin{aligned}
   dZ_t &=-\frac{1}{2}A^2Zdt+BdW_t, \\
   Z(0)&=0,\\
   \end{aligned}
   \right. 
  \end{equation}
  where $W$ is a $U$-cylindrical Wiener process and $U:=L_0^2(\mathbb{T}^2,\mathbb{R}^2)$.
  For $f\in L_0^2(\mathbb{T}^2,\mathbb{R}^2)$ we denote its component functions by $f_1,f_2\in L_0^2(\mathbb{T}^2)$ i.e. $f(x)=(f_1(x),f_2(x)),\forall x\in \mathbb{T}^2$. There exist two independent $L^2(\mathbb{T}^2)$-cylindrical Wiener processes $W^1$ and $W^2$ such that $W=(W^1,W^2)$. Set
  \begin{equation}\label{3.2.1}
  D(B)=H^1(\mathbb{T}^2,\mathbb{R}^2),B=\textrm{div},D(B^{\ast})=H^1_2 (\mathbb{T}^2),
  B^{\ast}=-\nabla.
  \end{equation}
   We know that
  $$Z_t(x)=\int_0^t e^{-\frac{t-s}{2}A^2}BdW_s=\int_0^t\langle K(t-s,x-\cdot),dW_s\rangle_U,$$
  where $K(t,x)=-\nabla_xM(t,x)=(K^1,K^2)$, and $M(t,x)$ is the kernel of $e^{-\frac{t}{2}A^2}$, that is, $M(t,x)=\sum_{k}e^{-\frac{t}{2}\lambda_k^2}e_k(x)$. 
  
  For any function $f$ on $\mathbb{T}^2$ , we can view it as a periodic function on $\mathbb{R}^2$ by defining
  $\bar{f}(x):=f(x+m)$, when $x+m\in\mathbb{T}^2$, $x\in \mathbb{R}^2$, $m\in\mathbb{Z}^2$.
  Moreover, define
  $$
\bar{K}^j(t,x):=-\mathcal{F}^{-1}(\pi i\xi_je^{-\frac{t}{2}|\pi\xi|^4})(x), j=1,2.  
  $$
   and $\bar{K}:=(\bar{K}^1,\bar{K}^2)$. By the Poisson summation formula (see \cite[Section \Rmnum{7}.2]{Stein:1972jf}) we know that
  \begin{equation}\label{3.3.01}
  K(t,x)=\sum_m\bar{K}(t,x+m), \forall t
  \end{equation}
and for any $f\in L^2(\mathbb{T}^2)$, $j=1,2$, $x\in \mathbb{T}^2$
\begin{equation}\label{3.4.01}
\begin{aligned}     
\partial_je^{-\frac{t}{2}A^2}f(x)&=\int_{\mathbb{T}^2}K^j(t,x-y)f(y)dy\\
&=\int_{\mathbb{R}^2}K^j(t,x-y)f(y)\mathbbm{1}_{\mathbb{T}^2}(y)dy\\
 &=\sum_m\int_{\mathbb{R}^2}\bar{K}^j(t,x-y+m)f(y)\mathbbm{1}_{\mathbb{T}^2}(y)dy\\
&=\int_{\mathbb{R}^2}\bar{K}^j(t,x-y)\sum_m\mathbbm{1}_{\mathbb{T}^2}(y+m)f(y+m)dy\\
 &=(\bar{K}^j(t,\cdot)*\bar{f})(x)  
\end{aligned},
\end{equation}
where we used (\ref{3.3.01}) in the third inequality and $\mathbbm{1}_{\mathbb{T}^2}$ is the indicator function of $\mathbb{T}^2$.
  Since 
  $$
 \bar{K}^j(t,x)=-\mathcal{F}^{-1}(\pi i\xi_je^{-\frac{t}{2}|\pi\xi|^4})(x)=t^{-\frac{3}{4}}\bar{K}^j(1,t^{-\frac{1}{4}}x)
  $$
  and 
  $$
 |\bar{K}^j(1,t^{-\frac{1}{4}}x)|\lesssim|\mathcal{F}^{-1}(\pi i\xi_je^{-\frac{1}{2}|\pi\xi|^4})(t^{-\frac{1}{4}}x)|\lesssim|1+t^{-\frac{1}{4}}x|^{-3},  
  $$
  we have the following estimate:
  \begin{equation}\label{3.3.1}
  |\bar{K}(t,x)|\lesssim t^{-\frac{\varepsilon}{4}}|x|^{-3+\varepsilon}, \forall \varepsilon\in [0,3].
  \end{equation}

   \begin{lemma}\label{l3.1.1}
   $Z\in C([0,T];\mathcal{C}^{-\alpha})$ $\mathbb{P}$-almost-surely, for all $\alpha>0$.
   \end{lemma}
   \proof By the factorization method in \cite{DaPrato:2004fs} we have that for $\kappa\in(0,1)$
$$Z(t)=\frac{\sin(\pi\kappa)}{\pi}\int_0^t(t-s)^{\kappa-1}\langle M(t-s,x-\cdot),U(s)\rangle ds,$$
and
$$U(s,\cdot)=\int_0^s(s-r)^{-\kappa} e^{-\f {s-r} 2 A^2}BdW_r.$$
A similar argument as in the proof of Lemma 2.7 in \cite{DaPrato:2004fs} implies that it suffices to prove that for $p>1/(2\kappa),$ 
\begin{equation}\label{5.2}
\mathbb{E}\|U\|_{L^{2p}(0,T;\mathcal{C}^{-\alpha})}<\infty.
\end{equation}
In fact, by (\ref{2.4}) we have that
$$\aligned \mathbb{E}\|U(s)\|_{-\alpha}^{2p}\lesssim 
&\sum_{\psi\in{\Psi}_\star}\sum_{n\geq0}\sum_{x\in\Lambda_n}\mathbb{E}2^{-2\alpha pn+2np}
|\langle U(s)
,\psi_{x}^{n}\rangle|^{2p}w(x)^{2p}\\
\lesssim&\sum_{\psi\in{\Psi}_\star}\sum_{n\geq0}\sum_{x\in\Lambda_n}2^{-2\alpha pn+2np}
(\mathbb{E}|\langle U(s)
,\psi_{x}^{n}\rangle|^{2})^pw(x)^{2p}.\endaligned$$
Here  we used that $\langle U(s)
,\psi_{x}^{n}\rangle$ belongs to the first order Wiener-chaos as well as Gaussian hypercontractivity (cf. \cite[Section 1.4.3]{Nualart:1995er} and \cite{thefreemarkofffie:4RNWhXc2}) in the second inequality.
Moreover, we obtain that
$$
\aligned 
\mathbb{E}|\langle U(s),\psi_{x}^{n}\rangle|^{2}
=& \mathbb{E}|\langle U^1(s),\psi_{x}^{n}\rangle|^{2}+\mathbb{E}|\langle U^2(s),\psi_{x}^{n}\rangle|^{2}\\
\leq& \sum_{j=1}^2 \int\int|\psi_{x}^{n}(y)\psi_{x}^{n}(\bar{y})|\int_0^s(s-r)^{-2\kappa}\bar{K}^j * \bar{K}^j(s-r,y-\bar{y})drdyd\bar{y}\\
\lesssim&\int\int|\psi_{x}^{n}(y)\psi_{x}^{n}(\bar{y})|\int_0^s(s-r)^{-\frac{\varepsilon}{2}-2\kappa}|y-\bar{y}|^{-4+2\varepsilon}drdyd\bar{y}\\
\lesssim&2^{2n-2\varepsilon n}s^{1-2\kappa-\frac{\varepsilon}{2}},
\endaligned$$
where $$U^j(y)=\int_0^s(t-s)^{\kappa -1}\langle K^j(s-r,y-\cdot),dW^j_r\rangle,j=1,2$$ and in the second inequality we used (\ref{3.4.01}) and we also used \cite[Lemma 10.17]{Hairer:2014hd} and  (\ref{3.3.1}) to deduce that $|\bar{K}^j*\bar{K}^j(s-r,y-\bar{y})|\lesssim |s-r|^{-\frac{\varepsilon}{2}}|y-\bar{y}|^{-4+2\varepsilon}$. 

In fact, we can decompose $\bar{K}$ into $\bar{K}:=\bar{K}_{\delta}+\bar{K}_{\delta}^c$, where $\bar{K}_{\delta}$ is a compactly supported function and satisfies (\ref{3.3.1}), and $\bar{K}_{\delta}^{c}$ is a Schwartz function. Then $\bar{K}*\bar{K}=\bar{K}_{\delta}*\bar{K}_{\delta}+H$, where $H$ is a Schwartz function. By \cite[Lemma 10.17]{Hairer:2014hd} we have $\bar{K}_{\delta}*\bar{K}_{\delta}(t,x)\lesssim t^{-\frac{\varepsilon}{2}}|x|^{-4+2\varepsilon}$ and $\bar{K}*\bar{K}$ satisfies the same inequality.

Thus, we have
$$
\mathbb{E}\|U(s)\|_{-\alpha}^{2p}
\lesssim \sum_{n\geq0} 2^{(4-2\varepsilon-2\alpha)pn}s^{(1-2\kappa-\frac{\varepsilon}{2})p}.
$$
Let $\kappa$ be so small that $2-\alpha<\varepsilon<2-4\kappa+\frac{2}{p}$, which implies that 
$$
4-2\varepsilon-2\alpha<0, (1-2\kappa-\frac{\varepsilon}{2})p>-1.
$$
Then (\ref{5.2}) follows.
 
   $\hfill\Box$
\vskip.10in
  
  Note that $BB^{\ast}=-A$. Then by Fourier expansion it is easy to see that $Z_t\sim \mathcal{N}(0,Q_t)$, i.e. for any $h\in \mathcal{S}(\mathbb{T}^2)$
\begin{equation*}
\mathbb{E}e^{i{}_{\mathcal{S}}\!\langle h,Z_t\rangle _{\mathcal{S}'}}=\exp(-\frac{1}{2}\langle Q_th,h\rangle),
\end{equation*}
where $Q_t=(-A)^{-1}(I-e^{-\frac{t}{2}A^2})$.
 
 According to the definition of $V^{\alpha}$ and Lemma \ref{l2.1}, we have $\mathcal{C}^{-\alpha}\subset V^{-\alpha-\varepsilon}$ for any $\alpha, \varepsilon>0$. Then by Lemma \ref{l3.1.1}, $\mu_t$ is supported on $V_0^{-\alpha}$ for any $\alpha>0$ and letting $t\to \infty$, by \cite[Example, 3.8.13]{Bogachev:1998hm}, the law of $Z_t$ converges to the Gaussian measure $\mu= \mathcal{N}(0,Q)$, which is also supported on $V_0^{-\alpha}$.

In the following we are going to define the Wick powers both in the state space and the path space. 

Firstly,
we define the Wick powers on $L^2(\cal{S}'(\mathbb{T}^2),\mu)$.
\vskip.10in

\textbf{Wick powers on $L^2(\mathcal{S}'(\mathbb{T}^2),\mu)$}

$\mu$ is of course also a measure supported on $\mathcal{S}'(\mathbb{T}^2)$. We have the well-known (Wiener-It\^{o}) chaos decomposition $$L^2(\mathcal{S}'(\mathbb{T}^2),\mu)=\bigoplus_{n\geq0}\mathcal{H}_n,$$
where $\mathcal{H}_n$ is the Wiener chaos of order $n$ (cf. \cite[Section 1.1.1]{Nualart:1995er}).
Now we define the Wick powers by using approximations: for $\phi\in \mathcal{S}'(\mathbb{T}^2)$ define $$\phi_\varepsilon:=\rho_\varepsilon*\phi,$$ with $\rho_\varepsilon$ an approximate delta function on $\mathbb{R}^2$ given by
$$\rho_\varepsilon(x)=\varepsilon^{-2}\rho(\frac{x}{\varepsilon})\in \mathcal{D}, \int \rho=1.$$
Here the convolution means that we view $\phi$ as a periodic distribution in $\mathcal{S}'(\mathbb{R}^2)$.
 For every $n\in\mathbb{N}$ we set $$:\phi_\varepsilon^n::=c_\varepsilon^{n/2}P_n(c_\varepsilon^{-1/2}\phi_\varepsilon),$$
where $P_n,n=0,1,...,$ are the Hermite polynomials defined by the formula
$$P_n(x)=\sum_{j=0}^{[n/2]}(-1)^j\frac{n!}{(n-2j)!j!2^j}x^{n-2j},$$
and
$c_\varepsilon=\int\phi^2_\varepsilon\mu(d\phi)=\int\int G(z-y)\rho_\varepsilon(y)dy\rho_\varepsilon(z)dz$. Then $$:\phi_\varepsilon^n:\in \mathcal{H}_n.$$
Here and in the following $G$ is the Green function associated with $-A$ on $\mathbb{T}^2$. In fact by \cite[Section 6.1, Chapter \Rmnum{7}]{Stein:1972jf}, $$G(x)=\sum_{k\in \mathbb{Z}^2\setminus \{(0,0)\}}\f1{\lambda_k}e_k(x)\simeq -\log|x|,|x|\to 0,$$
and $G$ is continuously differentiable outside $\{0\}$.

For Hermite polynomial $P_n$ we have that for $s,t\in\mathbb{R}$
\begin{equation}\label{5.1}
P_n(s+t)=\sum_{m=0}^nC_n^mP_m(s)t^{n-m},
\end{equation}
where $C_n^m=\frac{n!}{m!(n-m)!}$.

 \vskip.10in
A direct calculation yields the following:
 \vskip.10in

\begin{lemma}\label{l5.1}
Let $\alpha>0$, $n\in\mathbb{N}$ and $p>1$. $:\phi_\varepsilon^n:$ converges to some element in $L^p(\mathcal{S}'(\mathbb{T}^2),\mu;\mathcal{C}^{-\alpha})$ as $\varepsilon\rightarrow0$. This limit is called the $n$-th Wick power of $\phi$ with respect to the covariance $Q$ and denoted by $:\phi^n:$.
\end{lemma}

\proof The proof is similar to that of \cite[Lemma 3.1]{Rockner:2015uh} since the Green function $G$ has the same regularity. Therefore we omit it here for simplicity.
$\hfill\Box$

 \vskip.10in
 
 \textbf{Wick powers on a fixed probability space}
 
 Now we fix a probability space $(\Omega,\mathcal{F},P)$ and consider a $U$-cylindrical Wiener process $W$. In the following we assume that $\mathcal{F}$ is the $\sigma$-field generated by $\{\langle W_t,h\rangle, h\in U,t\in\mathbb{R}^+\}$. We also have the well-known (Wiener-It\^{o}) chaos decomposition $$L^2(\Omega,\mathcal{F},P)=\bigoplus_{n\geq0}\mathcal{H}'_n,$$
where $\mathcal{H}'_n$ is the Wiener chaos of order $n$ (cf. \cite[Section 1.1.1]{Nualart:1995er}). We can define Wick powers of $Z(t)$ with respect to different covariances by approximations. Let 
   \begin{equation*}
		\begin{aligned}     
     Z_{\varepsilon}(t,x)=\rho_{\varepsilon}*Z_t
     &=\int_0^t \langle B^*e^{-\f {t-s} 2 A^2}\rho_{\varepsilon,x},dW_s\rangle_U\\
     &=\int_0^t\langle K_{\varepsilon}(t-s,x-\cdot),dW_s\rangle_U,
     \end{aligned},\
     \end{equation*}
where $\rho_{\varepsilon,x}=\rho_{\varepsilon}(x-\cdot)$, $K_{\varepsilon}(t,x)=(\rho_{\varepsilon}*K_t^1,\rho_{\varepsilon}*K_t^2)$ and
$$K_t^j=-\sum_{k} (i \pi k_j)e^{-\f t 2 \lambda_k^2}e_k,j=1,2.$$

For any $n\in\mathbb{N}$, we set
$$
: Z_{\varepsilon}^{n}(t) :_{Q_{t}} :=\left(c_{\varepsilon, t}\right)^{\frac{n}{2}} P_{n}\left(\left(c_{\varepsilon, t}\right)^{-\frac{1}{2}} Z_{\varepsilon}(t)\right) \in \mathcal{H}_{n}^{\prime},
$$
where $P_n, n = 0, 1, \cdots$, are the Hermite polynomials and
$c_{\varepsilon, t}=\left\|\mathbbm{1}_{[0, t]} {K}_{\varepsilon}\right\|_{L^{2}\left(\mathbb{R} \times \mathbb{T}^{2};\mathbb{R}^2\right)}^{2}$.

\vskip.10in

\begin{lemma}(\cite[Lemma 3.4]{Rockner:2015uh})\label{l5.4}
For  $\alpha>0$, $p>1$,  $n\in\mathbb{N}$, $:Z_\varepsilon^n:$ converges  in $L^p(\Omega,C([0,T];\mathcal{C}^{-\alpha}))$. The limit is called Wick power of $Z(t)$ of order $n$ with respect to the covariance $Q$ and is denoted by $:Z^n(t):$.
\end{lemma}

\proof
The kernel $K$ is a little different from the kernel in \cite{Rockner:2015uh}. But (\ref{3.3.1}) satisfies the condition in \cite[Lemma 3.2]{Rockner:2015uh} and \cite[Lemma 4.1]{Zhu:2015uf} which leads to a similar proof as for \cite[Lemma 3.3]{Rockner:2015uh}, so we omit it here.
$\hfill\Box$

\begin{remark}
Here we do not combine the initial value with the Wick powers as in \cite{Mourrat:2015uo,Rockner:2015uh}, since we can obtain existence of solutions to the shifted equation (\ref{6.5}) for any initial value in $V_0^{-1}$ (see Section \ref{s4}).
\end{remark}
\vskip.10in

\textbf{Relations between two different Wick powers}

We introduce the following probability measure. Set $:q(\phi):=\frac{1}{4}:\phi^4:$,   $:p(\phi):=:\phi^{3}:$. Let $$\nu=c\exp(-N)\mu,$$
where $c$ is a normalization constant and $N={}_{\mathcal{S}'}\langle :q:,e_0\rangle_{\mathcal{S}}$. Then according to \cite[Lemma \Rmnum{5}.5 and Theorem \Rmnum{5}.7] {Simon:1974tu}  we have for every $p\in [1,\infty)$, $\varphi(\phi):=e^{-N}\in L^p(\mathcal{S}'(\mathbb{T}^2),\mu)$. 
\vskip.10in
The following result is about the relation between the two different Wick powers.
\vskip.10in

\begin{lemma} \label{l5.7}
Let $\phi$ be a measurable map from $(\Omega, \mathcal{F},\mathbb{P})$ to $C([0,T],B^{-\gamma}_{2,2})$ with $\gamma>2$, $\mathbb{P}\circ \phi(t)^{-1}=\nu$ for every $t\in[0,T]$ and let ${Z}(t)$ be defined as above. Assume in addition that  $y=\phi-{Z}\in  C((0,T];\mathcal{C}^{\beta})$ $\mathbb{P}$-a.s. for some $\beta>\alpha>0$. Here $C((0,T];\mathcal{C}^{\beta})$ is equipped with the norm $\sup_{t\in[0,T]}t^{\frac{\beta+\alpha}{4}}||\cdot||_{\beta}$. Then for every $t>0$, $n\in\mathbb{N}$
\begin{equation}\label{rwp}
:\phi^n(t):=\sum_{k=0}^nC_n^ky^{n-k}(t):{Z}^{k}(t):\quad \mathbb{P}-a.s..
\end{equation}
Here the Wick power on the left hand side is the limit obtained and defined in  Lemma \ref{l5.1}.
\end{lemma}

\proof By Lemma \ref{l5.4} it follows that for every $k\in \mathbb{N}$, $p>1$
$$:{Z}_\varepsilon^k:\rightarrow :{Z}^k:\quad \textrm{ in } L^p(\Omega, C((0,T];\mathcal{C}^{-\alpha})), \textrm{ as }\varepsilon\rightarrow0.$$
Since $y_\varepsilon=\phi_\varepsilon-{Z}_{\varepsilon}=\rho_\varepsilon*y$ and $y\in  C((0,T];\mathcal{C}^{\beta})$ $\mathbb{P}$-a.s., it is obvious that  $y_\varepsilon\rightarrow y$ in $C((0,T];\mathcal{C}^{\beta-\kappa})$ $\mathbb{P}$-a.s. for every $\kappa>0$ with $\beta-\kappa-\alpha>0$, which combined with Lemma \ref{l2.3} implies that for $k\in\mathbb{N}$, $k\leq n$,
$$y_\varepsilon^{n-k}:{Z}_\varepsilon^{k}:\rightarrow^{\mathbb{P}} y^{n-k}:{Z}^{k}: \quad \textrm{in } C((0,T];\mathcal{C}^{-\alpha}),  \textrm{ as }\varepsilon\rightarrow0.$$
Here $\rightarrow^{\mathbb{P}}$ means convergence in probability.
Since $e^{-N}\in L^p(\mathcal{S}'(\mathbb{T}^2),\mu)$ for every $p\geq 1$, by H\"{o}lder's inequality and Lemma \ref{l5.1} we get that for $t>0$ and $p>1$
$$:\phi_\varepsilon^n(t):\rightarrow :\phi^n(t): \quad \textrm{ in } L^p(\Omega,\mathcal{C}^{-\alpha}),  \textrm{ as }\varepsilon\rightarrow0.$$
Moreover, by (\ref{5.1}) we have
$$\aligned &:\phi_\varepsilon^n:=:(y_\varepsilon+{Z}_\varepsilon)^n:=c_\varepsilon^{n/2}P_n(c_\varepsilon^{-1/2}(y_\varepsilon+{Z}_\varepsilon))
\\=&\sum_{k=0}^nC_n^k c_\varepsilon^{n/2}P_k(c_\varepsilon^{-1/2}{Z}_\varepsilon)(c_\varepsilon^{-1/2}y_\varepsilon)^{n-k}\\=&\sum_{k=0}^nC_n^k:{Z}_\varepsilon^k: y_\varepsilon^{n-k},
\endaligned$$
which implies the result by letting $\varepsilon\rightarrow0$.
$\hfill\Box$
\vskip.10in

\section{The Solution to the Shifted Equation}\label{s4}

Now we fix a stochastic basis $(\Omega,\mathcal{F}, \{\mathcal{F}_t\}_{t\in[0,\infty)},\mathbb{P})$ and on it a $U$-cylindrical Wiener process $W$. Define ${Z}(t)=\int_0^te^{-(t-s)A^2/2}BdW(s)$ as in Section \ref{s3}. Now we consider the following shifted equation:
\begin{equation}\label{6.5}
     \left\{
   \begin{aligned}
   \frac{dY}{dt} &=-\frac{1}{2}A^2Y+\frac{1}{2}A\sum_{k=0}^{3}C_3^kY^{3-k}:{Z}^k:, \\
   Y(0)&=x .\\
   \end{aligned}
   \right.
  \end{equation}
  
Generally we consider initial data $x$ that are $\mathcal{F}_0$ measurable and belong to $V_0^{-1}, \mathbb{P}-a.s.$. To prove the existence of the solution to equation (\ref{6.5}), we use a smooth approximation on each path:
\begin{equation}\label{6.6}
   \left\{
   \begin{aligned}
   \frac{dY_{\varepsilon}}{dt} &=-\frac{1}{2}A^2Y_{\varepsilon}+\frac{1}{2}A\sum_{k=0}^{3}C_3^kY_{\varepsilon}^{3-k}:{Z}_{\varepsilon}^k:, \\
   Y_{\varepsilon}(0)&=x_{\varepsilon}, \\
   \end{aligned}
   \right.
  \end{equation}
  where ${Z}_{\varepsilon}={Z}*\rho_{\varepsilon}$, $x_{\varepsilon}=x*\rho_{\varepsilon}$, and $\rho_{\varepsilon}$ is as introduced in Section \ref{s3}. Note that the solution $Y$ to equation (\ref{6.5}) and the solution $Y_{\varepsilon}$ to (\ref{6.6}) also satisfy:
\begin{equation}\label{4.3.1}
\frac{dm(Y(t))}{dt}=0, m(Y(0))=0 ,
\end{equation}
  which means that $m(Y(t))=m(Y_{\varepsilon}(t))=0$.
  
  From Lemma \ref{l5.1} we know that there exists a $\Omega'\subset\Omega$, $\mathbb{P}(\Omega')=1$, such that for any $\omega\in\Omega'$, $Z(\omega), :Z^n:(\omega)\in C((0,T];\mathcal{C}^{-\alpha}), n=2,3$, $\forall \alpha>0$. Since ${Z}_{\varepsilon}(\omega)$ is smooth, by monotonicity trick in \cite[Theorem 5.2.2 and Theorem 5.2.4]{Liu:2015vb}, there exists a unique solution $Y_{\varepsilon}(\omega)$ to equation (\ref{6.6}) in $L^2(0,T;V^2_0)\cap C([0,T];L^2_0)$ for each $\omega\in \Omega'$. We are going to find a convergent subsequence of $\{Y_{\varepsilon}(\omega)\}$, which converge to a solution to equation (\ref{6.5}) and prove uniqueness of solutions to (\ref{6.5}). Then we obtain a unique $\mathcal{F}_t$-adapted solution to equation (\ref{6.5}).
  
   In this section we never distinguish $V^{\alpha}$, $H_2^{\alpha}$ and $B_{2,2}^{\alpha}$ since they have equivalent norms. For convenience we denote all of them by $H^{\alpha}$.

\begin{theorem}[a-priori estimate]\label{t6.5}
If $Y$ is a solution to equation (\ref{6.5}), then there exists a constant $C_T$ which only depends on $T$ and $Z(\omega)$, such that for $\forall t\in [0,T]$
\begin{equation}\label{4.4}
||Y||_{H^{-1}}^2-||x||_{H^{-1}}+\frac{1}{2}\int_0^t\left( ||Y(s)||_{H^1}^2+||Y(s)||_{L^4}^4\right) ds\leq C_T.
\end{equation}
Moreover there exist constants $C>0$, $\gamma_k>0$, $k=1,2,3$, such that for every $t\in(0,T]$
\begin{equation}\label{4.14.2}
\|Y_t\|_{H^{-1}}^2\leq C\left(t^{-1}\vee\left(\sum_{k=1}^3t^{-\rho\gamma_k}\sup_{0\leq r\leq t}\left(r^{\rho\gamma_k}\|:Z_r:\|_{-\alpha}^{\gamma_k}\right)\right)^{\frac{1}{2}}\right),
\end{equation}
where $\rho>0$ is the small enough constant introduced in Lemma \ref{l5.4}.
\end{theorem}

\proof Since
$$
\frac{dY}{dt}=-\frac{1}{2}A(AY-\sum_{k=0}^{3}C_3^kY^{3-k}:{Z}^k:),
$$
and $m(Y)=0$, taking scalar product with $(-A)^{-1}Y$ we obtain that
\begin{equation*}
\frac{d}{dt}||Y||_{H^{-1}}^2+||Y||_{H^1}^2+||Y||_{L^4}^4=-\langle\sum_{k=1}^{3}C_3^kY^{3-k}:{Z}^k:,Y\rangle,
\end{equation*}
that is
\begin{equation}\label{6.7}
\frac{d}{dt}||Y||_{H^{-1}}^2+||Y||_{H^1}^2+||Y||_{L^4}^4\lesssim |\langle Y,:{Z}^3:\rangle|+|\langle Y^2,:{Z}^2:\rangle|+|\langle Y^3,{Z}\rangle|.
\end{equation}
So we only need to estimate the right hand side of (\ref{6.7}). We only consider $|\langle Y^3,{Z}\rangle|$. The other terms can be estimated similarly. Lemma \ref{l2.3} implies
$$
|\langle Y^3,{Z}\rangle|\lesssim ||{Z}||_{-\alpha}||Y^3||_{B_{1,1}^{\alpha}},\quad \forall \alpha>0.
$$
Moreover, by Lemma \ref{l2.1} and Lemma \ref{l2.4}. 
$$
||Y^3||_{B_{1,1}^{\alpha}}\lesssim ||\Lambda^{\beta_0}Y^3||_{L^{p_0}}\lesssim ||\Lambda^{\beta_0}Y^{\frac{3}{2}}||_{L^{p_1}}||Y^{\frac{3}{2}}||_{L^{q_1}},
$$
where $\beta_0>\alpha, p_0>1$ and $\frac{1}{p_0}=\frac{1}{p_1}+\frac{1}{q_1}$.

Choose $q_1\leq \frac{8}{3}$ and $p_1>\frac{8}{5}$. Then we have
$$
||Y^{\frac{3}{2}}||_{L^{q_1}}=||Y||_{L^{\frac{3}{2}q_1}}^{\frac{3}{2}}\lesssim ||Y||_{L^4}^{\frac{3}{2}}.
$$
For $||\Lambda^{\beta_0}Y^{\frac{3}{2}}||_{L^{p_1}}$, we have
$$
||\Lambda^{\beta_0}Y^{\frac{3}{2}}||_{L^{p_1}}
\lesssim ||\Lambda^{\beta_1}Y^{\frac{3}{2}}||_{L^{p_2}}
\lesssim ||\Lambda Y^{\frac{3}{2}}||_{L^{p_2}}^{\beta_1}||Y^{\frac{3}{2}}||_{L^{p_2}}^{1-\beta_1},
$$
where $1<p_2<p_1<2$, $\beta_0=\beta_1-2(\frac{1}{p_2}-\frac{1}{p_1})$, $\beta_1<1$ and we used Lemma \ref{l2.1} in the first inequality and Lemma \ref{l2.5} in the second inequality. For $||\Lambda Y^{\frac{3}{2}}||_{L^{p_2}}$, let $p_2<\frac{8}{5}$. Then we have
$$
||\Lambda Y^{\frac{3}{2}}||_{L^{p_2}}
\lesssim ||Y^{\frac{1}{2}}\nabla Y||_{L^{p_2}}
\lesssim ||Y||_{H^1}||Y^{\frac{1}{2}}||_{L^{\frac{2p_2}{2-p_2}}}
\lesssim ||Y||_{H^1}||Y||_{L^{\frac{p_2}{2-p_2}}}^{\frac{1}{2}}
\lesssim ||Y||_{H^1}||Y||_{L^4}^{\frac{1}{2}},
$$
where we used H\"older's inequality in the second inequality. Furthermore
$$
||Y^{\frac{3}{2}}||_{L^{p_2}}\lesssim ||Y||_{L^{\frac{3p_2}{2}}}^{\frac{3}{2}}\lesssim ||Y||_{L^4}^{\frac{3}{2}}.
$$

Combining the above estimates we get
$$
||Y^3||_{B_{1,1}^{\alpha}}\lesssim ||Y||_{L^4}^{3-\beta_1}||Y||_{H^1}^{\beta_1}.
$$
Combining this with Lemma \ref{l5.4}, we have
$$
|\langle Y^3,{Z}\rangle|\lesssim ||Y||_{L^4}^{3-\beta_1}||Y||_{H^1}^{\beta_1}t^{-\frac{\rho}{4}}\lesssim t^{-\frac{\rho}{4} \gamma}+\kappa\left(||Y||_{L^4}^4+||Y||_{H^1}^2 \right) ,
$$
where $\gamma=\frac{4}{1-\beta_1}$ and we used Young's inequality. Choosing $\rho$ to be so small that $\frac{\rho}{4} \gamma<1$, we can conclude that there exists $\gamma_k>0$, $k=1,2,3$ such that $\frac{\gamma_k\rho}{4}<1$ and
$$
\frac{d}{dt}||Y||_{H^{-1}}^2+\frac{1}{2}\left( ||Y||_{H^1}^2+||Y||_{L^4}^4\right) \lesssim \sum_{k=1}^3\|:Z^k:\|_{-\alpha}^{\gamma_k}\lesssim \sum_{k=1}^3t^{-\frac{\gamma_k\rho}{4}}.
$$
Hence (\ref{4.4}) follows. Moreover, since $\|Y\|_{H^{-1}}\lesssim \|Y\|_{L^4}$ we have that 
$$
\frac{d}{dt}||Y||_{H^{-1}}^2+\frac{1}{2}||Y||_{H^{-1}}^4\lesssim \sum_{k=1}^3\|:Z^k:\|_{-\alpha}^{\gamma_k}.
$$
By \cite[Lemma 3.8]{Tsatsoulis:2016wm}, we have
$$ 
\|Y_t\|_{H^{-1}}^2\lesssim t^{-1}\vee\left(\sum_{k=1}^3t^{-\rho\gamma_k}\sup_{0\leq r\leq t}\left(r^{\rho\gamma_k}\|:Z_r:\|_{-\alpha}^{\gamma_k}\right)\right)^{\frac{1}{2}}.
$$
$\hfill\Box$
\vskip.10in

Since the approximating equation (\ref{6.6}) obey the same a-prior estimate as (\ref{6.5}), by (\ref{4.4}) we deduce that the sequence $\{Y_{\varepsilon}\}$ is bounded in $L^{\infty}(0,T;H^{-1})\cap L^4([0,T]\times \mathbb{T}^2)\cap L^2(0,T;H^1)$. This implies that $\{AY_{\varepsilon}\}$ is bounded in $L^2(0,T;H^{-1})$ and $\{Y_{\varepsilon}^3\}$ is bounded in $L^{4/3}([0,T]\times \mathbb{T}^2)$. Moreover, Lemma \ref{l2.1} and Lemma \ref{l5.4} imply that $\{:{Z}_{\varepsilon}^3:\}$ is bounded in $L^p(0,T;H^{-\alpha})$ for any $\alpha>0,\varepsilon>0$ and $p>1$. Then we can prove the following lemma:

\begin{lemma}\label{l6.6}
$\{\frac{dY_{\varepsilon}}{dt}\}$ is bounded in $L^{p}(0,T;H^{-3})$, where $p\in (1,\frac{4}{3})$.
\end{lemma}

\proof According to the argument before, we only need to show that $\{Y_{\varepsilon}^2{Z}_{\varepsilon}\}$ and $\{Y_{\varepsilon}:{Z}_{\varepsilon}^2:\}$ are bounded in $L^{p}(0,T;H^{-1})$ when $p\in (1,\frac{3}{4})$.

We omit $\varepsilon$ if there is no confusion in this proof.

For $Y^2{Z}$ we have
$$
||Y^2{Z}||_{B_{2,\infty}^{-\alpha}}\lesssim ||Y^2||_{B_{2,\infty}^{\beta_0}}||{Z}||_{-\alpha}\lesssim ||Y^2||_{B_{2,1}^{\beta_0}}||{Z}||_{-\alpha},
$$
where $\beta_0>\alpha >0$, we used Lemma \ref{l2.5} in the first inequality and Lemma \ref{l2.1} in the second inequality. Furthermore,
$$
||Y^2||_{B_{2,1}^{\beta_0}}\lesssim ||\Lambda^{\beta_1}Y^2||_{L^2}\lesssim ||\Lambda^{\beta_1}Y||_{L^{p_0}}||Y||_{L^{q_0}},
$$
where $\beta_1>\beta_0$, $\frac{1}{p_0}+\frac{1}{q_0}=\frac{1}{2}$, we used Lemma \ref{l2.1} in the first inequality and Lemma \ref{l2.4} in the second inequality.
By Lemma \ref{l2.1},
$B_{q,2}^{s}\subset L^q$ for any $q\geq 1$ and $s>0$. Since $H^{\delta}\simeq B_{2,2}^{\delta}\subset B_{q,2}^{\delta-1+\frac{2}{q}}$ for $q\geq2$, the Besov interpolation in Lemma \ref{l2.4} implies that
\begin{equation}\label{6.8}
||Y||_{L^{q_0}}
\lesssim ||Y||_{B_{q_0,2}^s}
\lesssim ||Y||_{B_{q_0,2}^{\frac{2}{q_0}}}^{1-\frac{1}{q_0}+\frac{s}{2}}||Y||_{B_{q_0,2}^{\frac{2}{q_0}-2}}^{\frac{1}{q_0}-\frac{s}{2}}
\lesssim ||Y||_{H^1}^{1-\frac{1}{q_0}+\frac{s}{2}}||Y||_{H^{-1}}^{\frac{1}{q_0}-\frac{s}{2}}.
\end{equation}
For $||\Lambda^{\beta_1}Y||_{L^{p_0}}$, let $p_0\geq 2$. Then we use Lemma \ref{l2.1} and Sobolev interpolation to get
$$
||\Lambda^{\beta_1}Y||_{L^{p_0}}\lesssim ||Y||_{H^{\beta_2}}\lesssim ||Y||_{H^1}^{\frac{1+\beta_2}{2}}||Y||_{H^{-1}}^{\frac{1-\beta_2}{2}},
$$
where $\beta_1=\beta_2+\frac{2}{p_0}-1=\beta_2-\frac{2}{q_0}$. Thus we have
\begin{equation}\label{6.9}
||Y^2||_{B_{2,1}^{\beta_0}}
\lesssim ||Y||_{H^1}^{\frac{3}{2}+\frac{\beta_1}{2}+\frac{s}{2}}||Y||_{H^{-1}}^{\frac{1}{2}-\frac{\beta_1}{2}-\frac{s}{2}}.
\end{equation}
By Lemma \ref{l5.4} we deduce that
$$
||Y^2{Z}||_{B_{2,\infty}^{-\alpha}}
\lesssim ||Y||_{H^1}^{\frac{3}{2}+\frac{\beta_1}{2}+\frac{s}{2}}||Y||_{H^{-1}}^{\frac{1}{2}-\frac{\beta_1}{2}-\frac{s}{2}}t^{-\frac{\rho}{4}}.
$$
 For any $p\in(1,\frac{4}{3})$, let $\beta_1$ and $s$ be small enough such that $(\beta_1+s+3)p<4$. Then Young's inequality implies that there exists $\gamma>0$ such that
$$
||Y^2{Z}||_{B_{2,\infty}^{-\alpha}}^{p}\lesssim ||Y||_{H^1}^2+||Y||_{H^{-1}}^{\frac{4}{3}\gamma (\frac{1}{2}-\frac{\beta_1}{2}-\frac{s}{2})}t^{-\frac{\rho}{3} \gamma}.
$$
For $\rho$ small enough, $\{Y_{\varepsilon}^2{Z}_{\varepsilon}\}$ is bounded in $L^{p}(0,T;B_{2,\infty}^{-\alpha})$.

On the other hand, 
$$
||Y:{Z}^2:||_{B_{2,\infty}^{-\alpha}}\lesssim ||Y||_{B_{2,\infty}^1}||:{Z}^2:||_{-\alpha}\lesssim ||Y||_{H^1}t^{-\frac{\rho}{4}},
$$
where we used Lemma \ref{l2.5} in the first inequality and Lemma \ref{l2.1}, Lemma \ref{l5.4} in the second inequality. Then by Young's inequality
$$
||Y:{Z}^2:||_{B_{2,\infty}^{-\alpha}}^{\frac{4}{3}}\lesssim ||Y||_{H^1}^2+t^{-\rho}.
$$
Choosing $\rho$ small enough we deduce that $\{Y_{\varepsilon}:{Z}_{\varepsilon}^2:\}$ is bounded in $L^{\frac{4}{3}}(0,T;B_{2,\infty}^{-\alpha})$.
By Lemma \ref{l2.1} we have $B_{2,\infty}^{-\alpha}\subset H_2^{-\alpha-\delta}$ for any $\delta>0$. Hence $\{Y_{\varepsilon}^2{Z}_{\varepsilon}\}$ and $\{Y_{\varepsilon}:{Z}_{\varepsilon}^2:\}$ are bounded in $L^{p}(0,T;H^{-1})$, $\forall p\in (1,\frac{4}{3})$, which implies the results.
$\hfill\Box$

\begin{theorem}\label{t6.7}
 For every $x\in V^{-1}_0$, there exists at least one solution to equation (\ref{6.5}) in $C([0,T];V^{-1}_0)\cap L^4([0,T]\times\mathbb{T}^2)\cap L^2(0,T;V^1_0)$.
\end{theorem}

\proof Since $H^1\subset H^{\delta}$ compactly for any $\delta<1$ (see \cite[Proposition 4.6]{Triebel:2006vl}), a classical compactness argument (cf. \cite[Lemma C.2]{Goldys:2009fb} or \cite[Theorem 2.1, Chapter \Rmnum{3}]{Temam:2001wj}) implies that there exists a sequence $\{\varepsilon_k\}$ and $Y\in L^{\infty}(0,T,H^{-1})\cap L^2(0,T;H^1)\cap L^4([0,T]\times \mathbb{T}^2)$, such that $Y_{\varepsilon_k}\to Y$ in $L^2(0,T;H^{\delta})\cap C([0,T];H^{-3})$, $\forall \delta<1$.

It is sufficient to show that for a suitable $\delta\in (0,1)$, the limit $Y$ we obtained above is a solution in $H^{-3}$. 

In fact, if $Y$ is a solution in $H^{-3}$, i.e. for any $h\in H^{3}$
\begin{equation}\label{6.10}
{}_{H^{-3}}\!\langle Y_t-Y_0,h\rangle_{H^3}
=-\frac{1}{2}\int_0^t {}_{H^{-1}}\!\langle A^2h,Y_s\rangle_{H^1}ds
+\frac{1}{2}\int_0^t {}_{H^{-1}}\!\langle \sum_{k=0}^{3}C_3^kY_s^{3-k}:{Z}_s^k:,Ah\rangle_{H^1}ds,
\end{equation}
$Y$ is in $L^{\infty}(0,T,H^{-1})\cap L^2(0,T;H^1)\cap L^4([0,T]\times \mathbb{T}^2)$. Then we take the scalar product of $\frac{dY}{dt}$ and $(-A)^{-1}Y$, which is just the duality in $H^{-3}$ and $H^3$. Hence
\begin{equation*}
\frac{d}{dt}||Y||_{H^{-1}}^2+||Y||_{H^1}^2+||Y||_{L^4}^4=-\langle\sum_{k=1}^{3}C_3^kY^{3-k}:{Z}^k:,Y\rangle.
\end{equation*}
Thus $||Y||_{H^{-1}}$ is continuous w.r.t $t$. Moreover, \cite[Lemma 1.4, Chapter \Rmnum{3}]{Temam:2001wj} implies that $Y$ is weakly continuous in $H^{-1}$. Hence $Y\in C([0,T];H^{-1})$.

We still write $\varepsilon$ instead of $\varepsilon_k$ if there is no confusion. Since $Y_{\varepsilon}$ is a solution to equation (\ref{6.6}), letting $\varepsilon\to 0$,
it is easy to see that
\begin{align*}
\lim_{\varepsilon\to 0}{}_{H^{-3}}\!\langle Y_{\varepsilon},h\rangle_{H^3}&= {}_{H^{-3}}\!\langle Y,h\rangle_{H^3},\quad \lim_{\varepsilon\to 0}{}_{H^{-1}}\!\langle A^2h,Y_{\varepsilon}\rangle_{H^1}= {}_{H^{-1}}\!\langle A^2h,Y\rangle_{H^1},\\
\lim_{\varepsilon\to 0} {}_{H^{-1}}\!\langle:\bar{Z}^3_{\varepsilon}:,Ah\rangle_{H^1}&= {}_{H^{-1}}\!\langle:\bar{Z}^3:,Ah\rangle_{H^1}.
\end{align*}

It remains to show that for any $h\in H^{1}$
\begin{align}
\lim_{\varepsilon\to 0}|\int_0^t \langle Y_{\varepsilon}^3(s)-Y^3(s),h\rangle ds|&=0, \label{6.11}\\
\lim_{\varepsilon\to 0}|\int_0^t  \langle Y_{\varepsilon}^2(s){Z}_{\varepsilon}(s)-Y^2(s){Z}(s),h\rangle ds|&=0, \label{6.12}\\
\lim_{\varepsilon\to 0}|\int_0^t \langle Y_{\varepsilon}(s):{Z}_{\varepsilon}^2:(s)-Y(s):{Z}^2:(s),h\rangle ds|&=0. \label{6.13}
\end{align}
 
Since $Y_{\varepsilon}\to Y$ in $L^4([0,T]\times\mathbb{T}^2)$, which is equivalent to $||Y_{\varepsilon}||_{L^4([0,T]\times\mathbb{T}^2)}\to ||Y||_{L^4([0,T]\times\mathbb{T}^2)}$ and $Y_{\varepsilon}\Rightarrow^m Y$, where $\Rightarrow^m$ means convergence in Lebesgue measure $m$ on $[0,T]\times\mathbb{T}^2$, we have $||Y_{\varepsilon}^3||_{L^{\frac{4}{3}}([0,T]\times\mathbb{T}^2)}\to ||Y^3||_{L^{\frac{4}{3}}([0,T]\times\mathbb{T}^2)}$ and $Y_{\varepsilon}^3\Rightarrow^m Y^3$. Then (\ref{6.11}) holds by uniform integrability.

For (\ref{6.12}), let $R_{\varepsilon}=Y_{\varepsilon}-Y$.  By the triangle inequality
$$
| \langle Y_{\varepsilon}^2{Z}_{\varepsilon}-Y^2{Z},h\rangle|
\lesssim |\langle R_{\varepsilon}(Y+Y_{\varepsilon})h,{Z}\rangle |+
|\langle {Z}_{\varepsilon}-{Z},Y^2h\rangle |.
$$

For the second term on the right hand side of the above inequality, we have
$$
|\langle {Z}_{\varepsilon}-{Z},Y^2h\rangle |
\lesssim ||{Z}_{\varepsilon}-{Z}||_{-\alpha}||Y^2h||_{B_{1,1}^{\alpha}}
\lesssim ||{Z}_{\varepsilon}-{Z}||_{-\alpha}||Y^2||_{B_{2,1}^{\alpha}}||h||_{B_{2,1}^{\alpha}},
$$
where we used Lemma \ref{l2.3} in the first inequality and Lemma \ref{l2.5} in the second inequality. By \cite[Remark 2, Section 3.2, Chapter 2]{Triebel:1992dx} we have $H^1\subset B_{2,1}^{\alpha}$ for any $\alpha<1$. Hence
$$
|\langle {Z}_{\varepsilon}-{Z},Y^2h\rangle |\lesssim ||{Z}_{\varepsilon}-{Z}||_{-\alpha}||Y^2||_{B_{2,1}^{\alpha}}||h||_{H^1}.
$$
Combining with (\ref{6.9}), we have
$$
|\langle {Z}_{\varepsilon}-{Z},Y^2h\rangle |\lesssim ||{Z}_{\varepsilon}-{Z}||_{-\alpha}||h||_{H^1}||Y||_{H^1}^{\frac{3}{2}+\frac{\beta_3}{2}+\frac{s}{2}}||Y||_{H^{-1}}^{\frac{1}{2}-\frac{\beta_3}{2}-\frac{s}{2}},
$$
where $\beta_3>\alpha>0, s>0$.
Let $\frac{3}{2}+\frac{\beta_3}{2}+\frac{s}{2}<2$. Then Lemma \ref{l5.4} and H\"older's inequality imply that
$$
|\int_0^t \langle {Z}_{\varepsilon}-{Z},Y^2h\rangle ds|\to 0,\quad \varepsilon \to 0.
$$
Similarly
$$
|\langle R_{\varepsilon}Yh,{Z}\rangle |\lesssim ||R_{\varepsilon}Y||_{B_{2,1}^{\alpha}}||h||_{H^1}||{Z}||_{-\alpha}.
$$
For $||R_{\varepsilon}Y||_{B_{2,1}^{\alpha}}$, we have
$$
||R_{\varepsilon}Y||_{B_{2,1}^{\alpha}}\lesssim ||R_{\varepsilon}Y||_{B_{2,2}^{\beta_0}}\lesssim ||\Lambda^{\beta_0} R_{\varepsilon}||_{L^4}||Y||_{L^4}+||\Lambda^{\beta_0}Y||_{L^4}||R_{\varepsilon}||_{L^4},
$$
where $\beta_0>\alpha>0$ and we used Lemma \ref{l2.1} in the first inequality and Lemma \ref{l2.4} in the second inequality. By Lemma \ref{l2.1} we have the Sobolev embedding $H_2^{\beta+\frac{1}{2}}\subset H^{\beta}_4$. Hence
 $$
 ||R_{\varepsilon}Y||_{B_{2,1}^{\alpha}}\lesssim ||R_{\varepsilon}||_{H^{\beta_0+\frac{1}{2}}}||Y||_{L^4}+||Y||_{H^{\beta_0+\frac{1}{2}}}||R_{\varepsilon}||_{H^{\frac{1}{2}}}.
 $$
 By Sobolev interpolation, choosing $\delta>\frac{1}{2}+\beta_0$, we have
 $$
||Y||_{H^{\beta_0+\frac{1}{2}}}\lesssim ||Y||_{H^1}^{\frac{3}{4}+\frac{\beta_0}{2}}||Y||_{H^{-1}}^{\frac{1}{4}-\frac{\beta_0}{2}}.
 $$
Moreover, since $\delta>\frac{1}{2}+\beta_0$, we have $||R_{\varepsilon}||_{H^{\frac{1}{2}}}\lesssim ||R_{\varepsilon}||_{H^{\delta}}$ and $||Y||_{H^{\frac{1}{2}+\beta_0}}\lesssim ||Y||_{H^{\delta}}$.
 Then we deduce that
 $$
  ||R_{\varepsilon}Y||_{B_{2,1}^{\alpha}}\lesssim ||R_{\varepsilon}||_{H^{\delta}}||Y||_{L^4}+||Y||_{H^1}^{\frac{3}{4}+\frac{\beta_0}{2}}||R_{\varepsilon}||_{H^{\delta}}||Y||_{H^{-1}}^{\frac{1}{4}-\frac{\beta_0}{2}}.
 $$
 Let $\beta_0<\frac{1}{2}$ such that
 $$
 \frac{3}{4}+\frac{\beta_0}{2}+1<2.
 $$
 Then by H\"older inequality, we get
 $$
 \int_0^t ||R_{\varepsilon}Y||_{B_{2,1}^{\alpha}}||h||_{H^1}||\bar{Z}||_{-\alpha} ds\lesssim \left(\int_0^t ||R_{\varepsilon}||_{H^{\delta}}^2 ds\right) ^{\frac{1}{2}}\left(\int_0^t (||Y||_{H^1}^2)Fds \right) ^{\frac{1}{2}}\left(\int_0^t ||Y||_{L^4}^4ds \right)^{\frac{1}{4}}\to 0 ,
 $$
 where $F\in L^{\infty}(0,T)$. 
 
Moreover, we have
 $$
|\langle Y_{\varepsilon}:{Z}_{\varepsilon}^2:-Y:{Z}^2:,h\rangle |
\lesssim  |\langle Y_{\varepsilon}(:{Z}_{\varepsilon}^2:-:{Z}^2:),h\rangle|
+|\langle R_{\varepsilon}:{Z}^2:,h\rangle|.
 $$
By essentially the same argument as above, (\ref{6.13}) also follows.
 
 Then we have got a solution $Y$ in $C([0,T];H^{-1})\cap L^4([0,T]\times\mathbb{T}^2)\cap L^2(0,T;H^1)$. Combining this with (\ref{4.3.1}), we have $Y\in C([0,T];V^{-1}_0)\cap L^4([0,T]\times\mathbb{T}^2)\cap L^2(0,T;V^1_0)$.
$\hfill\Box$
\vskip.10in

Now we have obtained the existence of solutions to equation (\ref{6.5}). The following is the uniqueness result.

\begin{theorem}\label{t6.8}
For every $x\in V^{-1}_0$,  there exists a unique solution to equation (\ref{6.5}) in $C([0,T];V^{-1}_0)\cap L^4([0,T]\times\mathbb{T}^2)\cap L^2(0,T;V^{1}_0)$.
\end{theorem}

\proof Suppose $u$, $v$ are two solutions of (\ref{6.5}) with the same initial value. Let $r=u-v$, then $r$ satisfies:
\begin{equation*}
  \left\{
   \begin{aligned}
   \frac{dr}{dt} &=-\frac{1}{2}A^2r+\frac{1}{2}A\sum_{k=0}^{3}C_3^k(u^{3-k}-v^{3-k}):{Z}^k: ,\\
   r(0)&=0 .\\
   \end{aligned}
   \right.
  \end{equation*}
Similarly to (\ref{6.7}) we have:
\begin{equation}\label{6.17}
\frac{d}{dt}||r||_{H^{-1}}^2+||r||_{H^1}^2\lesssim |\langle r^2(u+v),{Z}\rangle|+|\langle r^2,:{Z}^2:\rangle|.
\end{equation}

By Lemma \ref{l2.3} and Lemma \ref{l5.4} we know
$$
|\langle r^2,:{Z}^2:\rangle|\lesssim ||r^2||_{B_{1,1}^{\alpha}}t^{-\rho},
$$
where $\beta>\alpha>0$.
Then Lemma \ref{l2.1} and Lemma \ref{l2.4} imply that
$$
||r^2||_{B_{1,1}^{\alpha}}\lesssim ||\Lambda^{\beta_0}r^2||_{L^{\frac{4}{3}}}\lesssim ||\Lambda^{\beta_0}r||_{L^2}||r||_{L^4}\lesssim ||r||_{H^1}^{\frac{\beta_0+3}{2}}||r||_{H^{-1}}^{\frac{1-\beta_0}{2}},
$$
where $1>\beta_0>\alpha>0$ and we used the Sobolev interpolation and Sobolev embedding theorem in the last inequality. Then by Young's inequality, there exists a $\gamma_1>0$ such that for any $\varepsilon>0$
\begin{equation}\label{4.11.1}
|\langle r^2,:{Z}^2:\rangle|\lesssim \varepsilon||r||_{H^1}^2+||r||_{H^{-1}}^2t^{-\rho\gamma_1}.
\end{equation}
Let $\rho$ be small enough. Then $g:=t^{-\rho\gamma_1}\in L^1(0,T)$.

For $|\langle r^2(u+v),{Z}\rangle|$, we similarly obtain that
$$
|\langle r^2(u+v),{Z}\rangle|\lesssim ||r^2(u+v)||_{B_{1,1}^{\alpha}}||{Z}||_{-\alpha}\lesssim \left( ||ur^2||_{B_{1,1}^{\alpha}}+||vr^2||_{B_{1,1}^{\alpha}}\right) t^{-\rho}. 
$$
For $||ur^2||_{B_{1,1}^{\alpha}}$, we have
$$
||ur^2||_{B_{1,1}^{\alpha}}\lesssim ||\Lambda^{\beta_0}(ur^2)||_{L^{p_0}}\lesssim ||\Lambda^{\beta_0}u||_{L^{p_1}}||r^2||_{L^{q_1}}+||\Lambda^{\beta_0}r^2||_{L^{p_2}}||u||_{L^{q_2}}:=(\Rmnum{1})+(\Rmnum{2}),
$$
with $p_0>1$, $\beta_0>\alpha>0$, and $\frac{1}{p_0}=\frac{1}{p_1}+\frac{1}{q_1}=\frac{1}{p_2}+\frac{1}{q_2}$, $p_i,q_i>p_0, i=1,2$. Here we used Lemma \ref{l2.1} in the first inequality and Lemma \ref{l2.5} in the second inequality.

For $(\Rmnum{1})$,  according to (\ref{6.8}) we know that for any $s>0$
$$
||r^2||_{L^{q_1}}= ||r||_{L^{2q_1}}^2\lesssim ||r||_{H^1}^{2-\frac{1}{q_1}+s}||r||_{H^{-1}}^{\frac{1}{q_1}-s}.
$$
Moreover, let $p_1\geq 4$. Then
$$
||\Lambda^{\beta_0}u||_{L^{p_1}}\lesssim ||\Lambda^{\beta_1}u||_{L^4}\lesssim ||u||_{L^4}^{1-2\beta_1}||u||_{H^{\frac{1}{2}}_4}^{2\beta_1}\lesssim ||u||_{L^4}^{1-2\beta_1}||u||_{H^1}^{2\beta_1},
$$
where $\beta_1=\beta_0+\frac{1}{2}-\frac{2}{p_1}$ and we used Lemma \ref{l2.1} in the first inequality and Sobolev interpolation in the second inequality and Besov embedding Lemma \ref{l2.1} in the last inequality. 
Combining these estimates above we have
$$
(\Rmnum{1})\lesssim||r||_{H^1}^{2-\frac{1}{q_1}+s}||r||_{H^{-1}}^{\frac{1}{q_1}-s}||u||_{L^4}^{1-2\beta_1}||u||_{H^1}^{2\beta_1}.
$$
Hence by Young's inequality
$$
t^{-\rho}(\Rmnum{1})\lesssim \varepsilon ||r||_{H^1}^2+||r||_{H^{-1}}^2||u||_{H^1}^{\frac{4\beta_1}{\frac{1}{q_1}-s}}||u||_{L^4}^{\frac{2(1-2\beta_1)}{\frac{1}{q_1}-s}}t^{-\frac{2\rho}{\frac{1}{q_1}-s}}.
$$
Let $p_0$ be close\ to $1$ and $\beta_0$, $s$ be small enough such that $\frac{1}{p_1}>1-\frac{1}{p_0}+\beta_0+s$, which is equivalent to $\frac{2\beta_1}{\frac{1}{q_1}-s}+\frac{(1-2\beta_1)}{\frac{1}{q_1}-s}\frac{1}{2}< 1$. 
Then the H\"older inequality yields for $\rho$ small enough
$$
\int_0^t ||u||_{H^1}^{\frac{4\beta_1}{\frac{1}{q_1}-s}}||u||_{L^4}^{\frac{2(1-\beta_1)}{\frac{1}{q_1}-s}}{\tau}^{-\frac{2\rho}{\frac{1}{q_1}-s}}d\tau\lesssim \left( \int_0^t||u||_{H^1}^2d\tau\right)^{\frac{1}{2}} \left( \int_0^t||u||_{L^4}^4d\tau\right)^{\frac{1}{4}}.
$$
 Then we get 
$$f_1^u:=||u||_{H^1}^{\frac{4\beta_1}{\frac{1}{q_1}-s}}||u||_{L^4}^{\frac{2(1-\beta_1)}{\frac{1}{q_1}-s}}t^{-\frac{2\rho}{\frac{1}{q_1}-s}}\in L^1(0,T),$$
and for any $\varepsilon>0$,
\begin{equation}\label{6.18}
t^{-\rho}(\Rmnum{1})\lesssim \varepsilon ||r||_{H^1}^2+f_1^u ||r||_{H^{-1}}^2.
\end{equation}

For ($\Rmnum{2}$), let $q_2=4$. Then we have $\frac{1}{p_2}+\frac{1}{4}=\frac{1}{p_0}\in (\frac{3}{4}, 1)$, which implies that $p_2\in (\frac{4}{3},2)$.
Similarly by Lemma \ref{l2.5} we have
$$
||\Lambda^{\beta_0}r^2||_{L^{p_2}}\lesssim ||\Lambda^{\beta_0}r||_{L^{p_3}}||r||_{L^{q_3}},
$$
where $\frac{1}{p_3}+\frac{1}{q_3}=\frac{1}{p_2}$, $p_3,q_3>p_2$. From (\ref{6.8}) we know that for every $s>0$
$$
||r||_{L^{q_3}}
\lesssim ||r||_{H^1}^{1-\frac{1}{q_3}+\frac{s}{2}}||r||_{H^{-1}}^{\frac{1}{q_3}-\frac{s}{2}}.
$$
Let $p_3\geq 2$. Then by Lemma \ref{l2.1} we have
$$
||\Lambda^{\beta_0}r||_{L^{p_3}}\lesssim ||r||_{H^{\beta_2}}\lesssim ||r||_{H^1}^{\frac{1+\beta_2}{2}}||r||_{H^{-1}}^{\frac{1-\beta_2}{2}},
$$
where we used Sobolev interpolation in the second inequality and that $\beta_0=\beta_2-1+\frac{2}{p_3}$. Hence
$$
||\Lambda^{\beta_0}r^2||_{L^{p_2}}\lesssim ||r||_{H^1}^{\frac{3}{2}+\frac{\beta_2}{2}-\frac{1}{q_3}+\frac{s}{2}}||r||_{H^{-1}}^{\frac{1}{2}-\frac{\beta_2}{2}+\frac{1}{q_3}-\frac{s}{2}}
=||r||_{H^1}^{2+\frac{\beta_0}{2}-\frac{1}{p_2}+\frac{s}{2}}||r||_{H^{-1}}^{\frac{1}{p_2}-\frac{s}{2}-\frac{\beta_0}{2}}.
$$
Thus, we have
$$
(\Rmnum{2})\lesssim||r||_{H^1}^{2+\frac{\beta_0}{2}-\frac{1}{p_2}+\frac{s}{2}}||r||_{H^{-1}}^{\frac{1}{p_2}-\frac{s}{2}-\frac{\beta_0}{2}}||u||_{L^4}.
$$
Then by Young's inequality we have
$$
t^{-\rho}(\Rmnum{2})\lesssim \varepsilon||r||_{H^1}^2+||r||_{H^{-1}}^2||u||_{L^4}^{\frac{2}{\frac{1}{p_2}-\frac{s}{2}-\frac{\beta_0}{2}}}t^{-\frac{2\rho}{\frac{1}{p_2}-\frac{s}{2}-\frac{\beta_0}{2}}}.
$$
It is easy to see that $p_2<2$ yields $\frac{2}{\frac{1}{p_2}-\frac{s}{2}-\frac{\beta_0}{2}}\leq 4$ when $s$, $\beta_0$ are small enough. Then for small enough $\rho$ we have $f_2^u:=||u||_{L^4}^{\frac{2}{\frac{1}{p_2}-\frac{s}{2}-\frac{\beta_0}{2}}}t^{-\frac{2\rho}{\frac{1}{p_2}-\frac{s}{2}-\frac{\beta_0}{2}}}\in L^1(0,T)$.

Then we obtain that for any $\varepsilon>0$
$$
|\langle r^2u,Z\rangle|\lesssim \varepsilon||r||_{H^1}^2+f^u||r||_{H^{-1}}^2,
$$
where $f^u:=f_1^u+f_2^u\in L^1(0,T)$. 

The same holds with $u$ replaced by $v$. Let $f=f^u+f^v\in L^1(0,T)$. Then
$$
|\langle r^2(u+v),\bar{Z}\rangle|\lesssim \varepsilon||r||_{H^1}^2+f||r||_{H^{-1}}^2.
$$

Hence we get 
$$
\frac{d}{dt}||r||_{H^{-1}}^2+||r||_{H^1}^2\lesssim \varepsilon||r||_{H^1}^2+(f+g)||r||_{H^{-1}}^2.
$$
Choose a suitable $\varepsilon>0$ such that
\begin{equation}\label{uniqueness}
\frac{d}{dt}||r||_{H^{-1}}^2\lesssim (f+g)||r||_{H^{-1}}^2.
\end{equation}
Then by Gronwall's inequality we have
$$
||r(t)||_{H^{-1}}^2\lesssim ||r(0)||_{H^{-1}}^2\exp\left(\int_0^t f(s)+g(s)ds \right)=0.
$$
Since $V_0^{-1}$ is a subspace of $H^{-1}$, we obtain the uniqueness.

$\hfill\Box$
\vskip.10in

\begin{remark}\label{r4.7.2}
We emphasize that we can also obtain local well-posedness by using the fixed point argument in \cite{Debussche:2003du} and \cite{Mourrat:2015uo} with initial value in $\mathcal{C}^{-\frac{4}{3}+}$. Since we only have an $H^{-1}$-uniform estimate, the local solution cannot be extended to a global solution similarly as for the dynamical $\Phi_2^4$ equation.
\end{remark}

\section{Relation to the solution obtained by the Dirichlet form approach}\label{s5}

In Section \ref{s4}, we obtained a unique solution $Y$ to the shifted equation (\ref{6.5}). Now it is natural to ask whether $X:=Y+Z$ satisfies the original equation (\ref{1.1}) and having $\nu$ as an invariant measure. In this section, we
 are going to obtain a probabilistically weak solution of equation (\ref{1.1}) via the Dirichlet form approach and compare this solution with the solution we obtain in Section \ref{s4}. This helps us to obtain the uniqueness of the corresponding quasi-regular Dirichlet forms (see Theorem \ref{tmu}). As mentioned in the introduction, this may give us some hope to study the scaling limit of the Kawasaki dynamics of the Ising-Kac model.

First, we introduce the Gelfand triple that we will work on.
According to the definition of $V_0^{\alpha}$ and \cite[Theorem 3.1]{Hida:1980jt}, $\mu$ is supported on $V_0^{-s}$ for any $s>1$. So we fix a small enough $s_0>0$ and $V_0^{-1-s_0}$ as the state space and denote it by $E$ for convenience. By identifying $V_0^1$ and $V_0^{-1}$ via the Riesz isomorphism we have the following Gelfand triple:
\begin{equation}\label{5.1.1}
E^{\ast}\subset V_0^{-1}\subset E
\end{equation}
where $E^{\ast}=V^{s_0-1}_0$ and the dualization between $E^{\ast}$ and $E$ is ${}_{E^{\ast}}\langle u,v\rangle_{E}:={}_{V_0^{1+s_0}}\langle {Q}u,v\rangle_{V_0^{-1-s_0}}$ for any $u\in E^{\ast}, v\in E$. Here ${}_{V_0^{s}}\langle \cdot,\cdot \rangle_{V_0^{-s}}$ is given by
\begin{equation}\label{5.2.0}
{}_{V_0^{s}}\langle u,v\rangle_{V_0^{-s}}:=\sum_{k}{}_{\mathcal{S}'}\langle u,e_k\rangle_{\mathcal{S}}{}_{\mathcal{S}'}\langle v,e_k\rangle_{\mathcal{S}},u\in V_0^s, v\in V_0^{-s}.
\end{equation}
Then we have that
\begin{equation}\label{5.2.1}
{}_{E^{\ast}}\langle u,v\rangle_{E}=\langle u,v\rangle_{V_0^{-1}}, \forall u\in E^{\ast}, \forall v\in V_0^{-1}.
\end{equation}
Moreover, we define $\mathcal{F}C_b^{\infty};=\{f({}_{E^{\ast}}\langle l_1,\cdot \rangle_{E},\cdots,{}_{E^{\ast}}\langle l_m,\cdot\rangle _{E}):m\in \mathbb{N},f\in C_b^{\infty}(\mathbb{R}^m),l_1,\cdots,l_m\in E^{\ast}\}$. For all $\varphi =f({}_{E^{\ast}}\langle l_1,\cdot \rangle_{E},\cdots,{}_{E^{\ast}}\langle l_m,\cdot\rangle _{E})\in \mathcal{F}C_b^{\infty}$, we can define the directional derivative for $h\in V_0^{-1}$:
$$
\partial_h\varphi(z):=\lim_{t\to 0}\frac{\varphi(z+th)-\varphi(z)}{t}=\sum_{i=1}^m\partial_i f({}_{E^{\ast}}\langle l_1,\cdot \rangle_{E},\cdots,{}_{E^{\ast}}\langle l_m,\cdot\rangle _{E})\langle l_i,h\rangle_{V_0^{-1}}.
$$
Then by the Riesz representation theorem, there exists a map $\nabla \varphi:E \to V_0^{-1}$ such that
$$\langle \nabla \varphi(z),h\rangle _{V_0^{-1}}=\partial_h\varphi(z),h\in V_0^{-1}.$$ 

\subsection{Solution given by Dirichlet forms}\label{s5.1}

Since ${Q}^{-1-s_0}:V_0^{1+s_0} \to V_0^{-1-s_0}$ is the Riesz isomorphism for $V_0^{1+s_0}$, i.e.
$${}_{V_0^{1+s_0}}\!\langle h,k\rangle _{V_0^{-1-s_0}}=\langle {Q}^{-1-s_0}h,k\rangle _{V_0^{-1-s_0}},$$
$\mu$ is in fact a Gaussian measure on Hilbert space $V_0^{-1-s_0}$, with covariance operator $C:={Q}^{2+s_0}$, that is 
$$\int_{V_0^{-1-s_0}}e^{i\langle h,z\rangle_{V_0^{-1-s_0}}}\mu(dz)=\langle Ch,h\rangle_{V_0^{-1-s_0}}.$$ Then we have the following integration by parts formula for $\mu$:

\begin{proposition}\label{p4.1}
For all $F\in\mathcal{F}C_b^{\infty}, h\in V_0^{3+s_0}$, we have
  \begin{equation}\label{4.9}
   \int \partial_{h}F d\mu=\int {}_{E^{\ast}}\langle A^2h,\phi\rangle_{E} F(\phi)\mu(d\phi).
  \end{equation}
\end{proposition}

\proof
First, by \cite[Section 1.2.4]{DaPrato:2002ev} we know the reproducing kernel of $(V_0^{-1-s_0},\mu)$ is $V_{\mu}:=C^{1/2}V_0^{-1-s_0}=V_0^1$. Then by \cite[Theorem 3.1, Chapter \Rmnum{2}]{Ma:1992ec} we have
  \begin{align*}
  \int \partial_{ h}Fd\mu 
  &=\int \langle C^{-1} h,\phi\rangle _{V_0^{-1-s_0}}F(\phi)\mu(d\phi)\\
  &=\int \langle {Q}^{-2-s_0}h,\phi\rangle _{V_0^{-1-s_0}}F(\phi)\mu(d\phi)\\
  &=-\int {}_{V_0^{1+s_0}}\langle Ah,\phi\rangle_{V_0^{-1-s_0}}F(\phi)\mu(d\phi)\\
  &=-\int {}_{E^{\ast}}\langle {Q}^{-1}Ah,\phi\rangle _{E}F(\phi)\mu(dz)\\
  &=\int {}_{E^{\ast}}\langle A^2h,\phi\rangle_{E} F(\phi)\mu(dz).
  \end{align*}
$\hfill\Box$
\vskip.10in

\begin{remark}\label{r5.2.0}
In fact, by a similar argument as in \cite[(9.1.32)]{Gllmm:1982bb}, (\ref{4.9}) still holds for $F\exp(-N)$, where $N={}_{\mathcal{S}'}\langle :q^n:,e_0\rangle_{\mathcal{S}}$, i.e. for all $F\in\mathcal{F}C_b^{\infty}, h\in V_0^{3+s_0}$
$$
 \int \partial_{ h}\left( F\exp(-N)\right) d\mu=\int {}_{E^{\ast}}\langle A^2h,\phi\rangle_{E} F(\phi)\exp(-N(\phi))\mu(d\phi)
$$
\end{remark}

Then for the Gibbs measure $\nu$ defined in Section \ref{s3}, we have the following integration by parts formula:
 
 \begin{proposition}\label{p6.1}
   For all $F\in \mathcal{F}C_b^{\infty}, h\in V_0^{3+s_0}$, we have
    \begin{equation}\label{6.2}
     \int \partial_{ h}F d\nu=\int \left({}_{E^{\ast}}\langle A^2h,\phi\rangle _{E}-{}_{E^{\ast}}\langle Ah,:\phi^3:\rangle _{E}\right)F(\phi)\nu(d\phi),
    \end{equation}
    where $:\phi^3:$ has been constructed in Lemma \ref{l5.1}.
  \end{proposition}
  
  \proof  
  Acoording to Proposition \ref{p4.1} and Remark \ref{r5.2.0}
  \begin{align*}
  \int\partial_{ h}Fd\nu&=c\int(\partial_{ h}F)\exp(-N)d\mu\\
  &=c\int[\partial_{ h}(F\exp(-N))+F\exp(-N)\partial_{ h}N]d\mu\\
  &=\int F(\phi)\left({}_{E^{\ast}}\langle A^2h,\phi\rangle _{E}-\partial_{ h}N(\phi) \right)\nu(d\phi) 
  \end{align*}
  
  By \cite[Theorem 4.1.1]{OBATA:2006vc},
  $$\partial_{ h}:{\phi}_{\varepsilon}^n(x):=n:{\phi}_{\varepsilon}^{n-1}(x):(\rho_{\varepsilon}*h)(x).$$
  Here $\partial_{ h}:{\phi}_{\varepsilon}^n(x):$ is defined as the directional derivative of the function $\phi\to :{\phi}_{\varepsilon}^n(x)$. Then 
  $$
\partial_{ h}N_{\varepsilon}(\phi)=\langle :\phi_{\varepsilon}^3:,h*\rho_{\varepsilon}\rangle, 
  $$
  where $N_{\varepsilon}(\phi):=\langle\frac{1}{4}:\phi_{\varepsilon}^4:,e_0\rangle$.
  Letting $\varepsilon \to 0$, due to the closablity of $\partial_{\Pi h}$ in $L^2(E,\mu)$,
  $$\partial_{ h}N(\phi)=\langle :\phi^3:,h\rangle=-{}_{E^{\ast}}\langle Ah,:\phi^3:\rangle _{E},$$
  which implies
  $$ \int \partial_{ h}F d\nu=\int \left({}_{E^{\ast}}\langle A^2h,\phi\rangle _{E}-{}_{E^{\ast}}\langle Ah,:\phi^3:\rangle _{E}\right)F(\phi)\nu(d\phi).$$
 $\hfill\Box$
\vskip.10in

  \begin{theorem}\label{t6.2}
   The bilinear form
     $$
      \mathcal{E}(\varphi,\psi) := \frac{1}{2}\int \langle \nabla \varphi,\nabla \psi\rangle_{V_0^{-1}}d\nu,
                            \forall \varphi,\psi\in \mathcal{F}C_b^{\infty},
     $$
     is closable on $L^2(E,\nu)$. Its closure is a symmetric quasi-regular Dirichlet form denoted by $\left(\mathcal{E},D(\mathcal{E})\right)$.
  \end{theorem}

\proof
 Let $h_k=\sqrt{\lambda_k}e_k$. Then $\{h_k\}_{k\in \mathbb{Z}^2\setminus \{(0,0)\}}$ is an orthonormal basis of $V_0^{-1}$. Then
  $$\mathcal{E}(\varphi,\psi)=\frac{1}{2}\sum_{k\in\mathbb{Z}^2\setminus \{(0,0)\} }\int  \partial_{h_k}\varphi \partial_{h_k}\psi d\nu,\forall \varphi,\psi\in \mathcal{F}C_b^{\infty},$$
  By Proposition \ref{p6.1} we have $\int \partial_{h_k}\varphi d\nu=-\int \varphi \beta_{h_k}d\nu$, where $\beta_{h_k}\in L^2(E,\nu)$ and
  $$\beta_{h_k}(\phi)=-{}_{E^{\ast}}\langle A^2h_k,\phi\rangle _{E}+{}_{E^{\ast}}\langle Ah_k,:\phi^3:\rangle _{E},\forall k\neq(0,0).$$
  According to \cite[Proposition 3.3, Chapter \Rmnum{2}]{Ma:1992ec}, any $\mathcal{E}_1$-Cauchy sequence $\{u_n\}\subset \mathcal{F}C_b^\infty$, i.e. $\lim_{n,m\to \infty}\mathcal{E}(u_n-u_m,u_n-u_m)+\|u_n-u_m\|_{L^2(E,\nu)}=0$,  has a unique limit in $L^2(E,\nu)$, which is also called  $\left(\mathcal{E},\mathcal{F}C_b^{\infty}\right)$ is closable on $L^2(E,\nu)$.
The closure $\left(\mathcal{E},D(\mathcal{E})\right)$ is a symmetric Dirichlet form. Moreover, by \cite[Proposition 4.2, Chapter \Rmnum{4}]{Ma:1992ec}, it is standard to prove that the capacity of $\left(\mathcal{E},D(\mathcal{E}) \right)$ is tight, and according to the fact that  $\mathcal{F}C_b^{\infty}$ is dense in  $L^2(E,\nu)$ and separates the points in  $L^2(E,\nu)$, this means that $\left(\mathcal{E},D(\mathcal{E})\right)$ is a quasi-regular Dirichlet form in the sense of Definition \ref{qrdf} or \cite[Definition 3.1]{Ma:1992ec}.
$\hfill\Box$
\vskip.10in

 \begin{theorem}\label{tdfmp}
  There is a conservative Markov diffusion process on another probability space $(\overline{\Omega},\overline{\mathcal{F}},(\overline{\mathbb{P}}^z)_{z\in E})$
   $$M=\left(\overline{\Omega},\overline{\mathcal{F}},\mathcal{M}_t,(\overline{X}(t))_{t \ge 0},(\overline{\mathbb{P}}^z)_{z\in E}\right),$$
 which is \emph{properly associated with}
   $\left(\mathcal{E},D(\mathcal{E})\right)$, i.e. for $u\in L^2(E,\nu)\cap\mathcal{B}_b(E)$, the transition semigroup $P_tu(z):=\overline{\mathbb{E}}^z[u(\overline{X}(t))]$ is $\mathcal{E}$-{quasi-continuous} for all $t>0$ and is a $\nu$-version of $T_tu$ where $T_t$ is the semigroup associated with $(\mathcal{E},D(\mathcal{E}))$. 
 \end{theorem}
   \proof
   Since $\left(\mathcal{E},D(\mathcal{E})\right)$ is a quasi-regular Dirichlet form on $L^2(E,\nu)$, it is a direct consequece of Theorem \ref{tqrdf} and \cite[Theorem 3.6]{Albeverio:1991hk}.
   
$\hfill\Box$
\vskip.10in

   In particular, by the construction of $M$ (see e.g. \cite[Chapter 7]{Fukushima:1994kz} or \cite[Section 3, Chapter \Rmnum{4}]{Ma:1992ec}), $M$ can be chosen as the canonical process on $E$.
Moreover, for the conservative Markov diffusion $M$ properly associated with $(\mathcal{E},D(\mathcal{E}))$, we say that a set $S\subset E$ is a \emph{properly $\mathcal{E}$-exceptional set} if $\nu(S)=0$ and
     $\overline{\mathbb{P}}^z\left(\overline{X}(t)\in E\setminus S,\forall t\ge 0\right)=1$ for $z\in E\setminus S$. For the relation between $\mathcal{E}$-exceptional set and properly exceptional set, we have the following theorem:
 \begin{theorem}(\cite[Theorem 4.1.1]{Fukushima:1994kz})\label{exception}
If $N$ is an $\mathcal{E}$-exceptional set, then $N$ is contained in a properly $\mathcal{E}$-exceptional set $S$. $S$ can be taken to be Borel measurable.
 \end{theorem}

As in \cite{Albeverio:1991hk}, we now derive the SPDE satisfied by $\overline{X}$ in the analytically weak form:

   \begin{theorem}\label{t6.3}
    There exists a map $\overline{W}:\overline{\Omega} \to C([0,\infty);C([0,\infty);V_0^{-1-s_0}(\mathbb{T}^2,\mathbb{R}^2)))$, and a {properly $\mathcal{E}$-exceptional set} $S\in\mathcal{B}(E)$, such that $\forall z\in E\setminus S$, 
     $\overline{W}$ is a $U$-cylindrical Wiener process on $(\overline{\Omega}, \mathcal{M}_t,\overline{\mathbb{P}}^{z})$ and the sample paths of the associated process
     $M=\left(\overline{\Omega},\overline{\mathcal{F}},\mathcal{M}_t,(\overline{X}(t))_{t \ge 0},(\overline{\mathbb{P}}^z)_{z\in E}\right)$ on $E$ satisfy the following: for $h\in V^{3+s_0}$,
    \begin{equation}\label{6.3}
      \begin{aligned}
        {}_{E^{\ast}}\langle h,{\overline{X}}(t)-{\overline{X}}(0)\rangle _{E}=& -\frac{1}{2}\int_0^t {}_{E^{\ast}}\langle A^2h,{\overline{X}}(s)\rangle _{E}ds\\
       &+\frac{1}{2}\int_0^t {}_{E^{\ast}}\langle Ah,:{\overline{X}}(s)^3:\rangle _{E}ds\\
      & +\int_0^t \langle B^{\ast}h,d\overline{W}_s\rangle _{V_0^{-1}(\mathbb{T}^2,\mathbb{R}^2)},\forall t\ge 0, {\overline{\mathbb{P}}}^z-a.s.,
      \end{aligned}
     \end{equation}
     where $B, B^{\ast}$ are defined as in (\ref{3.2.1}).
     Moreover, $\nu$ is an invariant measure for M in the sense that $\int P_tud\nu=\int ud\nu$ for $u\in L^2(E,\nu)\cap \mathcal{B}_b(E)$.
   \end{theorem}

\proof
    Let $u_h(\phi)={}_{E^{\ast}}\langle h,\phi\rangle _{E},h\in V_0^{3+s_0}$, and let $(\mathcal{L}, D(\mathcal{L}))$ be the generator of $\left( \mathcal{E},D(\mathcal{E})\right)$. For any $v\in D(\mathcal{E})$
\begin{equation*}    
    \begin{aligned}
    \frac{1}{2}\int \langle  \nabla u_h,\nabla v\rangle _{V_0^{-1}}d\nu
    &=-\frac{1}{2}\int \partial_{ h}v(\phi) \nu(d\phi)\\
    &=\frac{1}{2}\int \left({}_{E^{\ast}}\langle A^2h,\phi\rangle _{E}-{}_{E^{\ast}}\langle Ah,:\phi^3:\rangle _{E}\right)v(\phi)\nu(d\phi).
    \end{aligned}
    \end{equation*}
    Hence $v\mapsto \frac{1}{2}\int \langle  \nabla u_h,\nabla v\rangle _{V_0^{-1}}d\nu$ is continuous in $L^2(E,\nu)$. By definition of the generator of $\left( \mathcal{E},D(\mathcal{E})\right)$, $u_h\in D(\mathcal{L})$ and $\mathcal{L}u_h(\phi)=-\frac{1}{2}\left({}_{E^{\ast}}\langle A^2h,\phi\rangle _{E}-{}_{E^{\ast}}\langle Ah,:\phi^3:\rangle _{E}\right)$.

    By the well-known Fukushima's decomposition (see e.g. \cite[Theorem 4.3]{Albeverio:1991hk}),  we have for q.e. $z\in E$,
    $$u_h(\overline{X}_t)-u_h(\overline{X}_0)=M_t^h+\int_0^t \mathcal{L}u_h(\overline{X}_s)ds=M_t^h-\frac{1}{2}\int_0^t \left( {}_{E^{\ast}}\langle A^2h,\overline{X}_s\rangle _{E}-{}_{E^{\ast}}\langle Ah,:\overline{X}_s^3:\rangle _{E}\right) ds,$$
    where $M^h$ is an additive functional (Definition \ref{AF}) and also a martingale  with $\langle M^h\rangle _t=t\|h\|_{V_0^{-1}}^2$.
   By \cite[Proposition 4.5]{Albeverio:1991hk}, 
    $$\langle M^h\rangle _t=\int_0^t \langle  \nabla u_h(\overline{X}_s),\nabla u_h(\overline{X}_s)\rangle _{V_0^{-1}}ds=t\| h\|_{V_0^{-1}}^2.$$
     Now we identify $M^h$ as the component in direction $h$ of the conservative noise.
   For $f=B^{\ast}\bar{Q}h\in U$, with $h\in V_0^{-1}$, define $\overline{W}^f_t:=M_t^h$ and let $D:=span\{B^{\ast}{Q}e_k: k\in \mathbb{Z}^2\setminus\{(0,0)\}\}$. Since $\|B^{\ast}{Q}h\|_U^2=\| h\|_{V_0^{-1}}^2$, it is easy to check that $\langle \overline{W}^f,\overline{W}^g\rangle _t=t\langle f,g\rangle _U$ for $f,g\in D$, where $\langle \overline{W}^f,\overline{W}^g\rangle _t$ is the bracket process of $\overline{W}^f$ and $\overline{W}^g$. Moreover, $D$ is dense in $U$ and $\overline{W}_t^{\cdot}$ is $\mathbb{Q}$-linear on $D$, since the embedding $U\to V_0^{-1-s}(\mathbb{T}^2,\mathbb{R}^2)$ is Hilbert-Schmidt for any $s>0$. By \cite[Theorem 6.2]{Albeverio:1991hk}, there exist a map $\overline{W}:\overline{\Omega} \to C([0,\infty);V_0^{-1-s}(\mathbb{T}^2,\mathbb{R}^2))$, and a properly $\mathcal{E}$-exceptional set $S\in\mathcal{B}(E)$, i.e. $\nu(S)=0$ and
     $\overline{\mathbb{P}}^z\left(\overline{X}(t)\in E\setminus S,\forall t\ge 0\right)=1$ for $z\in E\setminus S$, such that $\forall z\in E\setminus S$,
     $\overline{W}$ is a $U$-cylindrical Wiener process on $(\overline{\Omega},\mathcal{M}_t,\overline{\mathbb{P}}^z)$ such that for any $f\in D$
     $$
     {}_{V_0^{-1-s}}\langle \overline{W},f\rangle_{V_0^{1+s}}=\overline{W}^f, \overline{\mathbb{P}}^z-a.s.,
     $$ 
     where ${}_{V_0^{-1-s}}\langle \cdot,\cdot\rangle_{V_0^{1+s}}$ is defined by (\ref{5.2.0}).
     In particular,
     \begin{equation*}
     \langle B^{\ast}h,\overline{W}_t\rangle_{V_0^{-1}(\mathbb{T}^2,\mathbb{R}^2)}
     =\langle \overline{W}_t,B^{\ast}{Q}h\rangle_U=M_t^h,    
     \end{equation*}
     and $\overline{W}=(\overline{W}^1,\overline{W}^2)$, where $\overline{W}^i:\overline{\Omega}\to C([0,\infty);E), i=1,2$ are two independent $L_0^2$-cylindrical Wiener processes under $\overline{\mathbb{P}}^z$ for any $z\in E\setminus S$.
     
     Moreover, it is easy to check that for the constant function $\mathbf{1}$ on $E$, i.e. $\mathbf{1}(z)=1,\;\forall z\in E$, 
   $$\mathbf{1}\in \mathcal{F}C_b^\infty,\quad \nabla \mathbf{1}\equiv 0,$$
 and
   $$
-\int \mathcal{L}\mathbf{1}vd\nu= \mathcal{E}(\mathbf{1},v)=0,\forall v\in D(\mathcal{E}).
   $$
   Thus we obtain that $\mathcal{L}\mathbf{1}=0$ and $T_t\mathbf{1}=1$ in $L^2(E,\nu)$. Then by the symmetry of the semigroup $T_t$,
   $$
\int P_tud\nu=\int T_tud\nu=\int u T_t\mathbf{1}d\nu=\int ud\nu,\;\;\forall u\in L^2(E,\nu)\cap\mathcal{B}_b(E).   
   $$
   This yields that $\nu$ is an invariant measure of $\overline{X}$.
$\hfill\Box$
\vskip.10in

\begin{remark}
We mention that the above Dirichlet form arguments can be easily extend to the infinite volume case. The Dirichlet form $(\Lambda,D(\Lambda))$ for equation (\ref{1.1}) on $\mathbb{R}^2$ can be directly constructed by the closure of the following bilinear form
\begin{equation}
\Lambda(\varphi,\psi)=\int\langle\nabla f,\nabla g\rangle_{\dot{H}^{-1}}d\nu,\forall\varphi,\psi\in\mathcal{F}C_b^\infty(\dot{H}^{-1-}),
\end{equation}
where $\dot{H}^s$ is the homogeneous Sobolev space of order $s$, and $\nabla$ is the gradient in $\dot{H}^{-1}$.  By a similar argument as before, it is easy to check $(\Lambda,D(\Lambda))$ is quasi-regular and one obtains again a probabilistically weak solution directly.
\end{remark}

\subsection{Relation between the two solutions}\label{s5.2}

 In the following we discuss the relation between ${M}$ constructed above and the solution of the shifted equation (\ref{1.3}). 
 For $\overline{W}$ constructed in Theorem \ref{t6.3} define
 ${\overline{Z}}(t):= \int_{0}^te^{-(t-s)A^2/2}Bd \overline{W}_s$. We will prove that the difference $\overline{Y}:=\overline{X}-{\overline{Z}}\in C([0,T];\mathcal{C}^{-\alpha})$ is a solution to equation (\ref{1.3}) by replacing $Z$ by $\overline{Z}$.
 Recall that in Section \ref{s4}, for every $W$, we constructed a corresponding strong solution $X:=Y+Z$. In particular for $\overline{W}$, we have a solution $\tilde{X}$. 
By the pathwise uniqueness of (\ref{1.3}) in  Theorem \ref{t6.8}, we prove $\tilde{X}=\overline{X}$. Thus the law of the solution constructed by the Dirichlet form approach is the same as the law of $X:=Y+Z$ given in Section \ref{s4}. This implies that $\nu$ is also an invariant measure of $X$.

We also mention that by Lemma \ref{l2.1}, $\mathcal{C}^{-\alpha}\subset V^{-1}\subset V^{-1-s_0}$ for $\alpha\in(0,1)$, $\mathcal{C}^{-\alpha}$ is Borel-measurable subset of $V^{-1-s_0}$ and $\nu(E\cap \mathcal{C}^{-\alpha})=1$.
\vskip.10in

\begin{theorem}\label{t6.4}
 Let $\alpha\in(0,\frac{1}{3})$, $\alpha<\beta<2-\alpha$.  There exists a properly $\mathcal{E}$-exceptional set $S_2\subset E$  such that for every  $z\in(\mathcal{C}^{-\alpha}\cap E)\setminus S_2$ under $\overline{\mathbb{P}}^z$,  $\overline{Y}:=\overline{X}-{\overline{Z}}\in C((0,T];\mathcal{C}^{\beta})\cap C([0,T];\mathcal{C}^{-\alpha}\cap V_0^{-1})$ is a solution to the following equation:
\begin{equation}\label{6.4}
\overline{Y}(t)=\frac{1}{2} \int_0^te^{-(t-s)A^2/2}A\sum_{l=0}^{3} C_{3}^l\overline{Y}(s)^l:{\overline{Z}}(s)^{3-l}:ds+e^{-\frac{t}{2}A^2}\overline{X}(0).
\end{equation}
Here $C((0,T];\mathcal{C}^{\beta})$ is equipped with the norm $\sup_{t\in[0,T]}t^{\frac{\beta+\alpha}{4}}||\cdot||_{\beta}$.
Moreover,
\begin{equation}\label{pr.ex}
\overline{\mathbb{P}}^z[\overline{X}(t)\in (\mathcal{C}^{-\alpha}\cap E)\setminus S_2, \forall t\geq0]=1 \textrm{ for } z\in (\mathcal{C}^{-\alpha}\cap E)\setminus S_2.
\end{equation}
\end{theorem} 
 
 \proof For $z\in E\setminus {S}$ under $\overline{\mathbb{P}}^z$
 we have that
$$\overline{X}(t)=\frac{1}{2}\int_0^te^{-(t-\tau)A^2/2}A:\overline{X}(\tau)^3:d\tau+{\overline{Z}}(t)+e^{-\frac{t}{2}A^2}\overline{X}(0),$$
where $S$ is the properly $\mathcal{E}$-exceptional set in Theorem \ref{t6.3}.
Since $\nu$ is an invariant measure for $\overline{X}$, by Lemma \ref{l2.1} and Lemma \ref{l5.1} we conclude that  for every $T\geq0$,  $p>1, \delta>0,$ with $ 2\delta-\alpha<0,$ and $p_0>1$ large enough
$$\aligned&\int \overline{\mathbb{E}}^z\int_0^T\|:\overline{X}(\tau)^3:\|_{{-\alpha}}^pd\tau\nu(dz)\lesssim\int \overline{\mathbb{E}}^z\int_0^T\|:\overline{X}(\tau)^3:\|_{B^{\delta-\alpha}_{p_0,p_0}}^pd\tau\nu(dz)
\\=& T\int\|:\phi^3:\|_{B^{\delta-\alpha}_{p_0,p_0}}^p\nu(d\phi)\lesssim T\int\|:\phi^3:\|_{{2\delta-\alpha}}^p\nu(d\phi)<\infty,\endaligned$$
which implies that  for  $\nu-a.s.\;z\in E\setminus S$,  $\overline{\mathbb{P}}^z$-a.s.
\begin{equation}\label{holder}
:\overline{X}(\cdot)^3:\in L^p(0,T;\mathcal{C}^{-\alpha}),\quad \overline{\mathbb{E}}^z\int_0^T\|:\overline{X}(\tau)^3:\|_{{-\alpha}}^pd\tau<\infty, \quad \forall p>1.
\end{equation}
Here we used Lemma \ref{l2.1} to deduce the first result.  The second, however, does not imply the first  directly because of (\ref{2.1}).
By the definition of Wick power, it is easy to check that $z\mapsto \overline{\mathbb{E}}^z\int_0^T\|:\overline{X}(\tau)^3:\|_{{-\alpha}}^pd\tau$ is quasi-continuous in the sense of Definition \ref{dqe} (see e.g. \cite[Chapter 4, Exercise 2.9]{Ma:1992ec}). By Definition \ref{dqe}, (\ref{holder}) holds $\overline{\mathbb{P}}^z$-a.s. for q.e.\;$z\in E$.
Then Theorem \ref{exception} implies that
 there exists a properly $\mathcal{E}$-exceptional set $S_1\supset{S}$ such that (\ref{holder}) holds $\overline{\mathbb{P}}^z$-a.s. for  $z\in E\setminus S_1$. 

Moreover, Lemma \ref{l2.2} implies that for $\alpha<\beta<2-\alpha$, $z\in (E\cap \mathcal{C}^{-\alpha})\setminus S_1$
$$\int_0^te^{-(t-\tau)A^2/2}A:\overline{X}(\tau)^3:d\tau\in C([0,\infty);\mathcal{C}^{\beta})\quad \overline{\mathbb{P}}^z-a.s..$$
Now by Lemma \ref{l2.2} we conclude that for $z\in (E\cap \mathcal{C}^{-\alpha})\setminus S_1$, $e^{-\frac{t}{2}A^2}\overline{X}(0)\in C([0,T],\mathcal{C}^{-\alpha})\cap C((0,T],\mathcal{C}^{\beta})$. Thus,
$$\overline{X}-{\overline{Z}}\in  C([0,T],\mathcal{C}^{-\alpha})\cap C((0,T],\mathcal{C}^{\beta})\quad \overline{\mathbb{P}}^z-a.s. .$$
Since $\overline{\mathbb{P}}^\nu\circ \overline{X}(t)^{-1}=\nu$, by Lemma \ref{l5.7} we conclude that under $\overline{\mathbb{P}}^\nu$, by Fubini's theorem $\overline{Y}:=\overline{X}-{\overline{Z}}$ satisfies (\ref{6.4}) and for $\nu$-a.e. $z\in E$ under $\overline{\mathbb{P}}^z$, $\overline{Y}:=\overline{X}-{\overline{Z}}$ satisfies (\ref{6.4}). Moreover, it is easy to check that $\mathcal{C}^{-\alpha}\subset V^{-1}$ and $\langle\overline{Y},e_0\rangle=0$. Then we obtain that $\overline{Y}\in  C([0,T],\mathcal{C}^{-\alpha}\cap V_0^{-1})\cap C((0,T],\mathcal{C}^{\beta}), \overline{\mathbb{P}}^z-a.s. $ for $\nu-a.e. \;z\in E\cap \mathcal{C}^{-\alpha}$.

In the following we prove that these results hold under $\overline{\mathbb{P}}^z$ for $z$ outside a properly $\mathcal{E}$-exceptional set. First we have ${\overline{Z}}\in C([0,\infty);\mathcal{C}^{-\alpha})$ $\overline{\mathbb{P}}^{\nu}$-a.s., which combined with $\overline{X}-{\overline{Z}}\in C([0,T],\mathcal{C}^{-\alpha})$ implies
 $$\overline{\mathbb{P}}^\nu[\overline{X}\in C([0,\infty),\mathcal{C}^{-\alpha})]=1.$$
We also have
$$\aligned\overline{{Y}}(s,t_0):=&\overline{X}(s+t_0)-\overline{Z}(s+t_0)=\frac{1}{2}\int_{t_0}^{t_0+s}e^{-(t_0+s-\tau)A^2/2}A:\overline{X}(\tau)^3:d\tau\\+&e^{-sA^2/2}(\overline{X}(t_0)-\overline{Z}(t_0))\in  C((0,\infty)^2;\mathcal{C}^{\beta})\quad \overline{\mathbb{P}}^\nu-a.s..\endaligned$$
 Similar arguments as in the proof of Lemma \ref{l5.7} imply that $\forall s>0, t_0\geq0$
 $$\aligned \overline{\mathbb{P}}^\nu(:\overline{X}(s+t_0)^3:=\sum_{l=0}^{3} C_{3}^l\bar{{Y}}(s,t_0)^l:\overline{Z}(s+t_0)^{3-l}:,\\ \overline{X}\in C([0,\infty),\mathcal{C}^{-\alpha}),\bar{{Y}}\in C((0,\infty)^2;\mathcal{C}^{\beta}))=1,\endaligned$$
In the following we use $I_{t,t_0}$ to denote the equality $$\aligned&\int_0^te^{-(t-s)A^2/2}A:\overline{X}(s+t_0)^3:ds\\=&\sum_{l=0}^{3} \int_0^t e^{-(t-s)A^2/2}A C_{3}^l\bar{{Y}}(s,t_0)^l:\overline{Z}(s+t_0)^{3-l}:ds.\endaligned$$
We say that $I_{t,t_0}$ holds if the above identity holds in $\mathcal{C}^\beta$ and the $C^\beta$-norms of both sides are finite.
 Then using Fubini's theorem we know that
  $$\aligned \overline{\mathbb{P}}^\nu(I_{t,t_0} \textrm{ holds }   \forall t\geq0, \forall t_0\in \mathbb{Q}^+, \overline{X}\in C([0,\infty);\mathcal{C}^{-\alpha}),\bar{{Y}}\in C((0,\infty)^2;\mathcal{C}^{\beta}))=1.\endaligned$$
  Here we used $\overline{X}\in C([0,\infty);\mathcal{C}^{-\alpha})$ for $\alpha<\frac{1}{3}$ to make the right hand side of $I_{t,t_0}$ meaningful.
 It is obvious that the right hand side of the first equality  is continuous in $C^\beta$ with respect to $t_0$.  Since $\int_0^te^{-(t-s)A^2/2}A:\overline{X}(s+t_0)^3:ds=\int_{t_0}^{t+t_0}e^{-(t-s+t_0)A^2/2}A:\overline{X}(s)^3:ds$ we know that $\int_0^te^{-(t-s)A^2/2}A:\overline{X}(s+t_0)^3:ds$ is also continuous with respect to $t_0$ and we obtain that 
 $$\aligned \overline{\mathbb{P}}^\nu(I_{t,t_0} \textrm{ holds }      \forall t,t_0\geq0, \overline{X}\in C([0,\infty);\mathcal{C}^{-\alpha}),{\overline{Y}}\in C((0,\infty)^2;\mathcal{C}^{\beta}))=1.\endaligned$$
 Then we can claim that  there exists a properly $\mathcal{E}$-exceptional set $S_2\supset S_1$ such that for  $z\in (\mathcal{C}^{-\alpha}\cap E)\setminus S_2$ under $\overline{\mathbb{P}}^z$
 $$\aligned \overline{\mathbb{P}}^z(\overline{X}\in C([0,\infty);\mathcal{C}^{-\alpha}))=1.\endaligned$$

 Indeed, define
 $$\aligned\overline{\Omega}_0:=&\{\omega:\overline{X}\in C([0,\infty);\mathcal{C}^{-\alpha}),:\overline{Z}^k:\in C([0,\infty);\mathcal{C}^{-\alpha}),k=1,2,3,I_{t,t_0} \textrm{ holds }   \forall t, t_0\in \mathbb{Q}^+ \},\endaligned$$
 which is measurable in $\overline{\Omega}$ which is taken as the canonical space on $E$.
 Let $\Theta_t:\overline{\Omega}\rightarrow\overline{\Omega}, t > 0$, be the canonical shift, i.e. $\overline{X}\circ \Theta_t=\Theta_t(\omega) = \omega(\cdot+t)=\overline{X}(\cdot+t),\omega\in \overline{\Omega}$.
 
 In the following we will show that there exists a properly $\mathcal{E}$-exceptional set $\tilde{S}_2\supset S_1$ such that for any $z\in E\setminus \tilde{S}_2$, $t\in \mathbb{Q}^+$
\begin{equation}\label{aaf}
 \mathbb{P}^z\left(\cap_{t\in\mathbb{Q}^+}\Theta_t^{-1}\overline{\Omega}_0\right)=\mathbb{P}^z\left(\overline{\Omega}_0\right).
\end{equation}
As seen in the proof of Theorem \ref{t6.3}, by the Fukushima's decomposition (see e.g. \cite[Theorem 4.3]{Albeverio:1991hk}), for any $k\in\mathbb{Z}^2$, $M^{e_k}=\langle \overline{W},e_k\rangle$ is an additive functional in the sense of Definition \ref{AF}. Thus we can find a properly $\mathcal{E}$-exceptional set $\tilde{S}_2\supset S_1$ and a set $\Lambda\subset\overline{\Omega}$, such that \(\mathbb{P}^{x}(\Lambda)=1\), \(\forall x \in E \backslash \tilde{S}_2, \Theta_{t} \Lambda \subset \Lambda, \forall t>0\), and moreover, for each \(\omega \in \Lambda\), \(M^{e_k}_{t+s}(\omega)=M^{e_k}_{s}(\omega)+M^{e_k}_{t}\left(\Theta_{s} \omega\right), \forall t, s \geq 0, \forall k\in\mathbb{Z}^2\), which implies that for any $\omega\in\Lambda$, $t,s\geq 0$,
$$
 \overline{W}_{t+s}(\omega)= \overline{W}_{s}(\omega)+ \overline{W}_{t}\left(\Theta_{s} \omega\right).
$$
Thus $\overline{Z}\circ \Theta_t(r):=\int_0^re^{-(r-s)A^2/2}Bd(\overline{W}\circ\Theta_t)_s=\overline{Z}_{t+r}-e^{-rA^2/2}\overline{Z}_t$.
 Then we can define the Wick power $:(\overline{Z}\circ \Theta_t)^k(r):$ by Hermite polynomials in a similar way as Lemma \ref{l5.4}. In particular, for any $\omega\in\Lambda$, $t,s,r,t_0\geq 0$, by direct calculation, we obtain that for any $k\in\mathbb{N}$
 \begin{equation}\label{shift}
 \begin{aligned}
(\overline{X}\circ \Theta_t)(r)&=\overline{X}(t+r),\\
:(\overline{Z}\circ \Theta_t)^k(r):&=:(\overline{Z}(t+r)-e^{-rA^2/2}\overline{Z}(t))^k:=\sum_{l=0}^kC_k^l:\overline{Z}^{k-l}(t+r):(-e^{-rA^2/2}\overline{Z}(t))^l,\\
(\overline{Y}\circ \Theta_t)(s,t_0)&=(\overline{X}\circ \Theta_t)(s+t_0)-(\overline{Z}\circ \Theta_t)(s+t_0)=\overline{Y}(s+t_0,t)+e^{-\frac{s+t_0}{2}A^2}\overline{Z}(t).
\end{aligned}
 \end{equation}
For any $\omega\in \overline{\Omega}_0\cap \Lambda$, the continuity of $\overline{X}\circ \Theta_t$ and $:(\overline{Z}\circ \Theta_t)^k(r):$ follows from the first two identities above. Moreover, since $I_{s,t_0}$ holds for $\omega$, we have
\begin{align*}
\sum_{l=0}^kC_k^l\overline{Y}(s+t_0,t)^{k-l}:\overline{Z}(s+t_0+t)^l:&=\sum_{l=0}^kC_k^l\left((\overline{Y}\circ \Theta_t)(s,t_0)-e^{-\frac{s+t_0}{2}A^2}\overline{Z}(t)\right)^{k-l}:\overline{Z}(s+t_0+t)^l:\\
&=\sum_{l=0}^kC_k^l(\overline{Y}\circ \Theta_t)(s,t_0)^{k-l}:(\overline{Z}\circ \Theta_t)(s+t_0)^l:,
\end{align*}
which implies that $I_{s,t_0}$ also holds for $\Theta_t\omega$ for any $t\in \mathbb{Q}^+$. Hence we have proved that
$$\Theta_t^{-1}\overline{\Omega}_0\cap \Lambda\supset \overline{\Omega}_0\cap\Lambda,\quad t\in \mathbb{Q}^+.$$

 On the other hand, since (\ref{shift}) holds for all $\omega\in \Lambda$, if $\omega\in\cap_{t\in \mathbb{Q}^+}\Theta_t^{-1}\overline{\Omega}_0 \cap \Lambda$, by the definition of $\overline{\Omega}_0$, we conclude that $\overline{X}(\omega)(t+\cdot),\overline{Z}(\omega)(t+\cdot)\in C([0,\infty);\mathcal{C}^{-\alpha})$ for any $t\in \mathbb{Q}^{+}$. Hence $\overline{X}(\omega),\overline{Z}(\omega)\in C((0,\infty);\mathcal{C}^{-\alpha})$. It is also easy to check that $\lim_{t\to 0}\overline{X}(t)=\overline{X}(0)$ and $\lim_{t\to 0}\overline{Z}(\omega)(t)=\overline{Z}(\omega)(0)=0$. Similarly we can show that $:\overline{Z}(\omega)^k:\in C([0,\infty);\mathcal{C}^{-\alpha})$. Thus we obtain that  $\left(\cap_{t\in \mathbb{Q}^+}\Theta_t^{-1}\overline{\Omega}_0\right)\cap\Lambda\subset \overline{\Omega}_0\cap\Lambda$. Hence $\left(\cap_{t\in \mathbb{Q}^+}\Theta_t^{-1}\overline{\Omega}_0\right)\cap\Lambda= \overline{\Omega}_0\cap\Lambda$.
 Since \(\mathbb{P}^{x}(\Lambda)=1\), \(\forall x \in E \backslash \tilde{S}_2\), (\ref{aaf}) follows.

It is clear that $\mathbb{P}^y(\overline{\Omega}_0)=1$ for $\nu-a.e.\;y\in E$ and $\mathbb{P}^{\cdot}(\overline{\Omega}_0)\in L^2(E,\nu)\cap \mathcal{B}_b(E)$.
By the Markov property and the conservativity of $\overline{X}$, we know that 
$$\overline{\mathbb{P}}^z(\Theta_t^{-1}\overline{\Omega}_0)=\mathbb{P}^z(\mathbb{P}^{\bar{X}_t}(\overline{\Omega}_0))=P_t(\mathbb{P}^{\cdot}(\overline{\Omega}_0))(z)=1,\;\text{for}\;\nu-a.e. z\in E$$
 which by Theorem \ref{tdfmp}, $z\mapsto \overline{\mathbb{P}}^z(\Theta_t^{-1}\overline{\Omega}_0)$ has an $\mathcal{E}$-quasi-continuous $\nu$-version on $E$. It follows that for every $t>0$
$$\overline{\mathbb{P}}^z(\Theta_t^{-1}\overline{\Omega}_0)=1\quad q.e. z\in E,$$
which by (\ref{aaf}) implies that
$$\overline{\mathbb{P}}^z(\overline{\Omega}_0)=1\quad q.e. z\in E.$$
By the same argument as before, and by Theorem \ref{exception},  there exists a properly $\mathcal{E}$-exceptional set $S_2\supset \tilde{S}_2$ such that outside $S_2$ the result holds. Moreover, since for $z\in (\mathcal{C}^{-\alpha}\cap E)\setminus S_2$ and $\overline{\mathbb{P}}^z-a.s.$, both sides of $I_{t,t_0}$ belong to $C((0,\infty)^2,\mathcal{C}^\beta)$ as a map on $\mathcal{C}^\beta$ w.r.t $(t,t_0)\in (0,\infty)^2$, we conclude that,  for $z\in (\mathcal{C}^{-\alpha}\cap E)\setminus S_2$
 $$\aligned \overline{\mathbb{P}}^z(\overline{X}\in C([0,\infty);\mathcal{C}^{-\alpha}),I_{t,t_0} \textrm{ holds }   \forall t, t_0\geq 0)=1.\endaligned$$

 Now $\overline{Y}$ satisfies (\ref{6.4}) $\overline{\mathbb{P}}^z$-a.s. for $z\in (\mathcal{C}^{-\alpha}\cap E)\backslash S_2$. Moreover, for $z\in (\mathcal{C}^{-\alpha}\cap E)\backslash S_2$, $\overline{Y}\in C([0,\infty);\mathcal{C}^{-\alpha})\cap C([0,T],\mathcal{C}^\beta), {\overline{Z}}\in C([0,\infty);\mathcal{C}^{-\alpha})$ $\overline{\mathbb{P}}^z$-a.s., which combined with the definition of properly $\mathcal{E}$-exceptional set,  implies that $$\overline{\mathbb{P}}^z[\overline{X}(t)\in (\mathcal{C}^{-\alpha}\cap E)\setminus S_2, \forall t\geq0]=1 \textrm{ for } z\in (\mathcal{C}^{-\alpha}\cap E)\setminus S_2.$$
$\hfill\Box$

\begin{corollary}\label{c5.7}
Let $X=Y+Z$ where $Y$ is the unique solution to (\ref{6.5}) on probability space $(\Omega,\mathcal{F},\mathbb{P})$. Then $\nu$ is an invariant measure of $X$.
\end{corollary}
\proof
By Theorem \ref{t6.4},  we know that for $z\in (\mathcal{C}^{-\alpha}\cap E)\setminus S_2$, under the measure $\overline{\mathbb{P}}^z$, $\overline{Y}:=\overline{X}-\overline{Z}$ satisfies the shift equation (\ref{6.4}) and has better regularity i.e. $\overline{Y}\in C([0,T];\mathcal{C}^{-\alpha}\cap V_0^{-1})$. By the pathwise uniqueness result in Theorem \ref{t6.8}, we know that  $\overline{\mathbb{P}}^z\circ\left(\overline{Y},\overline{Z}\right)^{-1}=\mathbb{P}\circ (Y,Z)^{-1}$, for $Y(0)=z\in(\mathcal{C}^{-\alpha}\cap E)\setminus S_2$. Thus  $\overline{\mathbb{P}}^z\circ\overline{X}^{-1}=\mathbb{P}\circ X^{-1}$, for $X(0)=z\in(\mathcal{C}^{-\alpha}\cap E)\setminus S_2$. Since $\nu$ is an invariant measure of $\overline{X}$ and $\nu(\mathcal{C}^{-\alpha}\cap E)=1$, $\nu$ is an invariant measure of $X$. 

$\hfill\Box$
\vskip.10in

\begin{remark}
If we have a probabilistically strong solution on infinite volume case. We can obtain the above results and $\nu$ is an invariant measure for the solution to equation (\ref{1.1}) on $\mathbb{R}^2$.
\end{remark}

\subsection{Markov uniqueness in the restricted sense}\label{s5.3}

In this subsection we will use the pathwise uniqueness results in Theorem \ref{t6.8} to prove Markov uniqueness in the restricted sense and the uniqueness of the martingale (probabilistically weak) solutions to (\ref{1.1}) if the solution has $\nu$ as an invariant measure.

 Let $\mathcal{E}^{\textrm{q.r.}}$  be the set of all quasi-regular Dirichlet forms $(\tilde{\mathcal{E}}, D(\tilde{\mathcal{E}}))$ (cf. [MR92]) on $L^2(E;\nu)$ such that $\mathcal{F}C_b^\infty\subset D(L(\tilde{\mathcal{E}}))$ and $\tilde{\mathcal{E}}=\mathcal{E}$ on $\mathcal{F}C_b^\infty\times \mathcal{F}C_b^\infty$.
Here for a Dirichlet form $(\tilde{\mathcal{E}}, D(\tilde{\mathcal{E}}))$ we denote its generator by $(L(\tilde{\mathcal{E}}), D(L(\tilde{\mathcal{E}})))$.

In the following we  consider the martingale problem in the sense of \cite{Albeverio:1994dn} and probabilistically weak solutions to (\ref{1.1}):
\vskip.10in

\begin{definition}(i) A continuous strong Markov process $M=(\Omega,\mathcal{F}, (\mathcal{M}_t), X_t, (\mathbb{P}^z))$ in the sense of [MR92, Chapter IV] with state space $E$ is said to solve the martingale problem for $(L(\mathcal{E}),  D)$ if for all $u\in  D$, $u(X(t))-u(X(0))-\int_0^t L(\mathcal{E})u(X(s))ds$, $t\geq0$, is an $(\mathcal{M}_t)$-martingale under $\mathbb{P}^\nu$.

 (ii) A continuous strong Markov process $M=(\Omega,\mathcal{F}, (\mathcal{M}_t), X_t, (\mathbb{P}^z))$  with state space $E$ is called a  probabilistically weak solution to (\ref{1.1}) if there exists two maps $W^i:\Omega\rightarrow C([0,\infty);E)$ $i=1,2$  such that for $\nu$-a.e. $z$ under $\mathbb{P}^z$, $W:=(W^1,W^2)$ is an $L^2({\mathbb{T}^2,\mathbb{R}^2})$- cylindrical Wiener process with respect to $(\mathcal{M}_t)$ and the sample paths of the associated  process satisfy (\ref{6.3})  for all $h\in V^{3+s_0}$.
\end{definition}

\begin{remark} 
If $M$ is a  probabilistically weak solution to (\ref{1.1}), we can easily check that it also solves the martingale problem. Conversely, if $M$ solves the martingale problem, then with the same arguement in Theorem \ref{t6.3}, there exists an $L_0^2({\mathbb{T}^2,\mathbb{R}^2})$-cylindrical Wiener process W such that $(X,W)$ satisfies (\ref{6.3})  for $h\in V^{3+s_0}$. That is, these two definitions are equivalent.
\end{remark}

To explain the uniqueness result below we also introduce the following  concept:

Two strong Markov processes $M$ and $M'$ with state space $E$ and transition semigroups $(p_t)_{t>0}$ and $(p_t')_{t>0}$ are called \emph{$\nu$-equivalent} if there exists $S\in\mathcal{B}(E)$ such that (i) $\nu(E\backslash S)=0$, (ii) $\mathbb{P}^z[X(t)\in S, \forall t\geq0]=\mathbb{P}'^z[X'(t)\in S, \forall t\geq0]=1,  z\in S$, (iii) $p_tf(z)=p_t'f(z)$ for all $f\in\mathcal{B}_b(E), t>0$ and $z\in S$.

\vskip.10in

Combining Theorem \ref{t6.8} and Theorem \ref{t6.4},  we  obtain Markov uniqueness in the restricted sense for $(L(\mathcal{E}), D)$ (see part (iii)) and the uniqueness of martingale (probabilistically weak) solutions to (1.1) having $\nu$ as an invariant measure (see part (i), (ii)):
\vskip.10in

\begin{theorem}\label{tmu}

 (i) There exists  (up to  $\nu$-equivalence) exactly one  probabilistically weak solution  $M$ to (\ref{1.1}) satisfying  $\mathbb{P}^z(X\in C([0,\infty);E))=1$ for $\nu$-a.e.  and having $\nu$ as an invariant measure, i.e. for the transition semigroup $(p_t)_{t\geq0}$, $\int p_t fd\nu=\int fd\nu$ for $f\in L^2(E;\nu)$.

(ii) $\sharp\mathcal{E}^{\textrm{q.r.}}=1$. Moreover,  there exists (up to $\nu$-equivalence) exactly one continuous strong Markov process $M$ with state space $E$  associated with a Dirichlet form $(\mathcal{E},D(\mathcal{E}))$ solving the martingale problem for $(L(\mathcal{E}),D)$.
\end{theorem}

\proof
The proof is the same as \cite[Theorem 3.12]{Rockner:2015uh}. We only explain the idea of proof here.

For (i), suppose the there is another probabilistically weak solution $\tilde{M}=(\tilde{\Omega},\tilde{\mathbb{P}}^z)$ to (\ref{1.1}). By the same argument as in the proof of Theorem \ref{t6.4}, $\tilde{X}-\int_0^{\cdot}e^{-(\cdot-s)A^2/2}Bd\tilde{W}_s$ also satisfies the shifted equation (\ref{6.4}) by replacing $\overline{Z}$ to $\tilde{Z}:=\int_0^{\cdot}e^{-(\cdot-s)A^2/2}Bd\tilde{W}_s$. Similarly to the proof of Corollary \ref{c5.7}, the pathwise uniqueness result in Theorem \ref{t6.8} implies that $\tilde{\mathbb{P}}^z\circ \tilde{X}^{-1}=\overline{\mathbb{P}}^z\circ \overline{X}^{-1}$ for $\nu-a.e. z\in E$. Thus
 $\tilde{M}$ is properly associated with the same quasi-regular Dirichlet form $(\mathcal{E},D(\mathcal{E}))$ as $\overline{M}$, then the assertion follows from \cite[Theorem 6.4]{Ma:1992ec}.

For (ii),  the second result follows from the first result and \cite[Theorem 3.4]{Albeverio:1994dn}. For the first result, suppose $\tilde{\mathcal{E}}\in \mathcal{E}^{\textrm{q.r.}}$ and there exists a unique Markov process $\tilde{M}$ associated with $\tilde{E}$. Similarly as before, the pathwise uniqueness for the shifted equation yields the uniqueness in law.
Thus the semigroup of $\tilde{M}$ are the same as $M$'s. Then $M$ and $\tilde{M}$ generate the same Dirichlet form, i.e. $(\mathcal{E},D(\mathcal{E})=(\tilde{\mathcal{E}},D(\tilde{\mathcal{E}}))$.
$\hfill\Box$

\begin{remark}\label{rss}
By the same argument, we can also prove the uniqueness of the probabilistically strong stationary solution to (\ref{1.1}). 
For more details we refer to \cite[Section 3.5]{Rockner:2015uh}.
\end{remark}

\section{Ergodicity}\label{s6}

Let $X=Y+Z$, where $Y$ is the solution to equation (\ref{6.5}). By the uniqueness of the solution $Y$ we have that $X$ is a Markov process. Let $P_t$ be the semigroup of $X$, i.e
$$
P_t\Phi(x)=\mathbb{E}\Phi\left(X(t,x)\right),\quad \forall \Phi\in C_b(V_0^{-1}).
$$

We recall that the $U$-cylindrical Wiener process $W$ takes values in $C([0,T], V_0^{-1-s_0}(\mathbb{T}^2,\mathbb{R}^2)), \mathbb{P}-a.s.$, for any $s_0>0$. Let $\mathcal{D}$ denote the Fr\'echet derivative of functions on $C([0,T], V_0^{-1-s_0}(\mathbb{T}^2,\mathbb{R}^2))$(i.e. with respect to the noise). We also denote the Cameron-Maritin space by $\mathcal{CM}:=\{\omega:\partial_t\omega\in L^2([0,T],L_0^2( \mathbb{T}^2;\mathbb{R}^2)),\omega(0)=(0,0)\}$. Here we view $\partial_t\omega$ as a function on $[0,T]\times \mathbb{T}^2$ rather than lying in the tagent space of $\mathbb{T}^2$.

\begin{proposition}\label{p62.1}
For a fixed $x\in V^{-1}_0$, let $\mathfrak{X}_t^x:=X(t,x)=Z_t+Y(t,x)$ be a map from $C([0,T], V_0^{-1-s_0})$ to $V^{-1}_0$. For any $\omega\in \mathcal{CM}$ its directional derivative $\mathcal{D}\mathfrak{X}_t^x(\omega)$ is given in mild form as
\begin{equation}
\mathcal{D}\mathfrak{X}_t^x(\omega)=\frac{1}{2}\int_0^te^{-(t-s)A^2/2}A\sum_{l=0}^2 3C^l_2 Y^{2-l}(s):Z^l_s:\mathcal{D}\mathfrak{X}_s^x(\omega)ds+\int_0^t e^{-(t-s)A^2/2}Bd\omega_s.
\end{equation}
\end{proposition}

The proof of Proposition \ref{p62.1} can be obtained by using approximation or the implicit function theorem (see \cite[Theorem 19.28]{Driver03}, \cite{Hairer:2016wl}, \cite{Tsatsoulis:2016wm})

Let $D$ denote the Fr\'echet derivative of functions on $V^{-1}_0$. We also consider the following equation:
\begin{equation}
\left\{
\begin{aligned}
\partial_tJ_{s,t}h&=-\frac{1}{2}A^2J_{s,t}h+\frac{1}{2}A\left(\sum_{l=0}^2 3C^l_2 Y^{2-l}(t):Z^l_t:J_{s,t}h\right)\\
J_{s,s}h&=h\in V^{-1}_0
\end{aligned}
\right..
\end{equation}
Then $J_{0,t}h=DX(t,x)(h)$, i.e. it is the derivative of $X(t,\cdot)$ in the direction $h$. For $\omega\in \mathcal{CM}$, by Duhamel's principle
\begin{equation}\label{62.3}
\mathcal{D}\mathfrak{X}_t^x(\omega)=\int_0^tJ_{s,t}B\partial_s\omega(s)ds.
\end{equation}

We define the stopping time
\begin{equation}
\tau^r:=\inf\{t\in(0,T): t^{\rho}||:Z^k_t:||_{-\alpha}>r,k=1,2,3\},
\end{equation}
where $\rho>0$ is the small enough constant introduced in Lemma \ref{l5.4}.

\begin{proposition}\label{p62.2}
For any  $x\in V^{-1}_0$ with $\|x\|_{H^{-1}}\leq R$, there exist  constants $C_1(R), C_2(R)$ such that for all $t\leq\tau^r$
$$
\sup_{s\leq t}||Y_s||_{H^{-1}}\vee\int_0^t||Y_s||_{L^4}^4ds\vee\int_0^t||Y_s||_{H^1}^2ds\leq C_1\quad \text{and}\quad \sup_{s\leq t}||J_{0,s}h||_{H^{-1}}\leq C_2||h||_{H^{-1}}
$$
\end{proposition}
\proof
The first bound with constant $C_1$ follows from the proof of Theorem \ref{t6.5}.

For the second bound, we note that $J_{0,t}h$ satisfies the following equation:
\begin{equation*}
\left\{
\begin{aligned}
\frac{du}{dt}&=-\frac{1}{2}A^2u+\frac{1}{2}A\left( \sum_{l=0}^2 3C^l_2 Y^{2-l}(t):Z^l_t:u\right) \\
u(0)&=h
\end{aligned}
\right..
\end{equation*}
Taking scalar product with $(-A)^{-1}u$, we obtain that
$$
\frac{d}{dt}||u||_{H^{-1}}^2+||u||_{H^1}^2=-3\langle Y^2+2YZ+:Z^2:,u^2\rangle,
$$
that is
$$
\frac{d}{dt}||u||_{H^{-1}}^2+||u||_{H^1}^2\leq 6|\langle YZ,u^2\rangle|+3|\langle :Z^2:,u^2\rangle|.
$$
Following the same argument that we used to get estimate (\ref{6.17}) and using the first bound, we use Gr\"onwall's inequality to obtain
the second bound.
$\hfill\Box$
\vskip.10in

Let $\chi_r \in C^{\infty}(\mathbb{R})$ such that $\chi_r(\zeta)\in[0,1]$ for all $\zeta\in\mathbb{R}$, and
\begin{equation*}
\chi_r(\zeta)=
\left\{
\begin{aligned}
1&, |\zeta|\leq\frac{r}{2}\\
0&, |\zeta|\geq r
\end{aligned}
\right..
\end{equation*}

Following the notation in \cite{Tsatsoulis:2016wm}, we set
\begin{equation}
C^{3,-\alpha}(0,T):=C([0,T];\mathcal{C}^{-\alpha})\times C((0,T];\mathcal{C}^{-\alpha})^2,
\end{equation}
and $\underline{Z}:=\left(Z,:Z^2:,:Z^3: \right) \in C^{3,-\alpha}(0,T)$. We also define
$$
\normmm{\underline{Z}}_t:=\max_{k=1,2,3}\left\lbrace \sup_{0\leq s\leq T}s^{\rho}||:Z_s^k:||_{-\alpha}\right\rbrace .
$$

\begin{theorem}(Bismut-Elworthy-Li Formula)\label{t62.3}
Let $x\in V^{-1}_0$, $\Phi\in C_b^1(V^{-1}_0)$ and $\omega$ be a process taking values in the Cameron-Martin space $\mathcal{CM}$ with $\partial_s\omega$ adapted. Assume that there exists a deterministic constant $C\equiv C(t)$ such that $||\partial_s\omega||_{L^2(0,t;U)}\leq C $ $\mathbb{P}-a.s.$. Then we have
\begin{equation}\label{62.6}
\begin{aligned}
\mathbb{E}[D\Phi(\mathfrak{X}_t^x)(\mathcal{D}\mathfrak{X}_t^x(\omega))\chi_r(\normmm{\underline{Z}}_t)]&=\mathbb{E}\left(\Phi(\mathfrak{X}_t^x)\chi_r(\normmm{\underline{Z}}_t)\int_0^t\partial_s\omega(s)\cdot dW_s \right) \\
&-\mathbb{E}\left(\Phi(\mathfrak{X}_t^x)\partial_{+}\chi_r(\normmm{\underline{Z}}_t)(\omega) \right) 
\end{aligned},
\end{equation}
where
\begin{equation}
\begin{aligned}
\partial_{+}\chi_r(\normmm{\underline{Z}}_t)(\omega)&=\partial_\zeta\chi_r(\normmm{\underline{Z}}_t)\partial_+\normmm{\underline{Z}}_t(\underline{Y})\\
\partial_+\normmm{\underline{Z}}_t(\underline{Y})&=\lim_{\delta\to 0^+}\frac{\normmm{\underline{Z}+\delta\underline{Y}}_t-\normmm{\underline{Z}}_t}{\delta}
\end{aligned},
\end{equation}
$\underline{Y}=\left(Q_{\omega}(\cdot),2ZQ_{\omega}(\cdot),3:Z^2:Q_{\omega}(\cdot) \right)\in C^{3,-\alpha}(0,t)$ and 
$$
Q_{\omega}(\cdot):=\int_0^{\cdot} e^{-(\cdot-s)A^2/2}B\partial_s\omega(s)ds.
$$

\end{theorem}
\proof 
This is proved by the same calculation as that in the proof of \cite[Theorem 5.4]{Tsatsoulis:2016wm}.

$\hfill\Box$
\vskip.10in

We use (\ref{62.6}) to prove the following proposition.

\begin{proposition}
There exists universal constants $\theta_1>0$ such that for every $T>0$, $x\in V_0^{-1}$ with $\|x\|_{H^{-1}}\leq R$, there exists a constant $C\equiv C(T,R)>0$ satisfying
\begin{equation}
|P_t\Phi(x)-P_t\Phi(y)|\leq C(T,R)\frac{1}{t^{\theta_1}}\|\Phi\|_{\infty}\|x-y\|_{H^{-1}}+2\|\Phi\|_{\infty}\mathbb{P}(t\geq\tau^{\frac{r}{2}})
\end{equation}
for every $y\in V^{-1}_0$, $\|x-y\|_{H^{-1}}\leq 1$, $\Phi\in C_b^1(V^{-1}_0)$ and $t\in [0,T]$.
\end{proposition}
\proof Let $\Phi\in C_b^1(V^{-1}_0)$. Then
$$
|P_t\Phi(x)-P_t\Phi(y)|=|E\left[\Phi\left( X(t,x)\right) - \Phi\left( X(t,y)\right)\right]|\leq I_1+I_2,
$$
where 
\begin{align*}
I_1&:=|\mathbb{E}\left[\Phi\left( X(t,x)\right) - \Phi\left( X(t,y)\right)\chi_r(\normmm{\underline{Z}}_t)\right]|\\
I_2&:=|\mathbb{E}\left[\Phi\left( X(t,x)\right) - \Phi\left( X(t,y)\right)\left( 1-\chi_r(\normmm{\underline{Z}}_t)\right) \right]|.
\end{align*}
For the second term we have that $I_2\leq 2\|\Phi\|_{\infty}\mathbb{P}(t\geq\tau^{\frac{r}{2}})$. By the mean value theorem we get that
\begin{align*}
I_1&=\Big|\mathbb{E}\left( \int_0^1D\Phi\left(\mathfrak{X}^{z_{\lambda}}_t(y-x)\right) d\lambda\cdot\chi_r(\normmm{\underline{Z}}_t) \right)\Big|\\
&=\Big|\int_0^1\mathbb{E}\left(D\Phi\left(\mathfrak{X}^{z_{\lambda}}_t \right) (y-x)\chi_r(\normmm{\underline{Z}}_t) \right)d\lambda \Big|,
\end{align*}
where $z_{\lambda}:=x+\lambda(y-x)$. For any $h\in V^{-1}_0$, let $\omega$ be such that $B\partial_s\omega(s)=J_{0,s}h$ for $s\leq\tau^r$ and $0$ otherwise. Then $\partial_s\omega(s)$ satisfies the condition in Theorem \ref{t62.3}. Furthermore, by (\ref{62.3}) and since $J_{0,s}J_{s,t}=J_{0,t}$, we have $\mathcal{D}\mathfrak{X}^{z_{\lambda}}_t(\omega)=tD\mathfrak{X}_t^{z_{\lambda}}(h)$. Then we use (\ref{62.6}) to obtain that
\begin{align*}
\mathbb{E}\left(D\left[\Phi(\mathfrak{X}_t^{z_{\lambda}}) \right](h)\chi_r(\normmm{\underline{Z}}_t)  \right)&=\frac{1}{t}\mathbb{E}\left(\Phi(\mathfrak{X}_t^{z_{\lambda}})\int_0^t\partial_s\omega(s)\cdot dW_s\chi_r(\normmm{\underline{Z}}_t) \right)  \\
&-\frac{1}{t}\mathbb{E}\left(\Phi(\mathfrak{X}_t^{z_{\lambda}})\partial_{+}\chi_r(\normmm{\underline{Z}}_t)(\omega) \right) .
\end{align*}
Then we have
\begin{align*}
I_1\leq \frac{1}{t}\|\Phi\|_{\infty}\int_0^1\mathbb{E}\Big|\int_0^t\partial_s\omega(s)dW_s\chi(\normmm{\underline{Z}}_t)\Big|d\lambda
+\frac{1}{t}\|\Phi\|_{\infty}\int_0^1\mathbb{E}\Big|\partial_+\chi_r(\normmm{\underline{Z}}_t)(\omega)\Big|d\lambda.
\end{align*}
For the first term we have
\begin{align*}
\mathbb{E}\Big|\int_0^t\partial_s\omega(s)dW_s\chi(\normmm{\underline{Z}}_t)\Big|&\leq \mathbb{E}\Big|\int_0^{t\wedge\tau^r}\partial_s\omega(s)\cdot dW_s\Big|\\
&\leq \left(\int_0^{t\wedge\tau^r}\|\partial_s\omega(s)\|_U^2ds \right)^{\frac{1}{2}} \\
&\lesssim \left(\int_0^{t\wedge\tau^r}\|J_{0,s}h\|_{H^{-1}}^2ds \right)^{\frac{1}{2}} \\
&\leq C_2t\|h\|_{H^{-1}},
\end{align*}
where we used the Cauchy-Schwarz inequality and It\^o's isometry in the second step and Propostion \ref{p62.2} in the last step.

By the definition of $\partial_+\chi_r(\normmm{\underline{Z}}_t)(\omega)$, we have
\begin{align*}
\Big|\partial_+\chi_r(\normmm{\underline{Z}}_t)(\omega)\Big|\leq\partial_+\normmm{\underline{Z}}_t(\underline{Y})
\leq \normmm{\underline{Y}}_t
\lesssim \normmm{\underline{Z}}_t\|Q_{\omega}(t)\|_{\beta},
\end{align*}
where $\underline{Y}$ is as introduced in Theorem \ref{t62.3} and we used Lemma \ref{l2.3} in the last inequality. Moreover, we use Lemma \ref{l2.2} and Lemma \ref{l2.5} to obtain
$$
\|Q_{\omega}(t)\|_{\beta}\lesssim\int_0^t(t-s)^{-\frac{\beta+2}{4}}\|J_{0,s}h\|_{-2}ds\lesssim\int_0^t(t-s)^{-\frac{\beta+2}{4}}\|J_{0,s}h\|_{H^{-1}}ds\lesssim C_2t^{\frac{2-\beta}{4}}\|h\|_{H^{-1}}.
$$
Choosing $\beta$ small enough, we deduce that there exists a constant $\theta_1\in(0,\frac{1}{2})$, such that
$$
I_1\lesssim C_2\frac{1}{t^{\theta_1}}\|\Phi\|_{\infty}\|h\|_{H^{-1}}.
$$
Letting $h=y-x$ we finish the proof.
$\hfill\Box$
\vskip.10in

We denote by $\|\mu_1-\mu_2\|_{TV}$ the total variation distance of two probability measures $\mu_1$, $\mu_2$ on $V^{-1}_0$ given by
$$
\|\mu_1-\mu_2\|_{TV}:=\sup_{\|\Phi\|_{\infty}\leq 1}\Big|\int\Phi d\mu_1-\int\Phi d\mu_2\Big|.
$$

\begin{theorem}\label{t6.5.2}
There exists $\theta\in (0,1)$ such that for any $x,y\in V_0^{-1}$ with $\|x\|_{H^{-1}}\leq R$ and $\|x-y\|_{H^{-1}}\leq 1$ there exists a constant $C(R)>0$ satisfying
$$
\|P_t(x,\cdot)-P_t(y,\cdot)\|_{TV}\leq C(R)\|x-y\|_{H^{-1}}^{\theta},
$$
for every $t\geq 1$.
\end{theorem}

\proof
The proof is the same as that of \cite[Theorem 5.8]{Tsatsoulis:2016wm}. Therefore we omit it.
$\hfill\Box$
\vskip.10in

In order to use Krylov-Bogoliubov method to prove the existence of an invariant measure, the $H^{-1}$ uniform estimate is not enough. We need to find a space compactly embedded in $H^{-1}$ where the solution is bounded in probability. We make use of the integrability on a smaller space, which is compactly embeded in $H^{-1}$. Thus we have:
\begin{theorem}\label{t6.62}
For every $x\in V^{-1}_0$, there exists a probability Borel measure $\nu_x$ on $V^{-1}_0$ such that $\nu_x$ is an invariant measure for the semigroup $\{P_t,t\geq 0\}$ on $V_0^{-1}$.
\end{theorem}
\proof
By (\ref{4.14.2}) and a similar argument as in the proof of \cite[Corollary 3.10]{Tsatsoulis:2016wm} we have that
\begin{equation}\label{4.15.2}
\sup_{x\in V_0^{-1}}\sup_{t>0}(t\wedge 1)\mathbb{E}\|X(t,x)\|_{H^{-1}}^2<\infty.
\end{equation}
By the uniqueness of the solution, we know $X(t,x)=Z_{t-1,t}+Y_{t-1,t}$, where $Z_{s,t}:=\int_s^te^{-(t-r)A^2/2}BdW_r$ and $Y_{s,r}$, $r\geq t-1$, solves the equation
\begin{equation}\label{6.2.2}
     \left\{
   \begin{aligned}
   \frac{dY_{s,r}}{dr} &=-\frac{1}{2}A^2Y_{s,r}+\frac{1}{2}A\sum_{k=0}^{3}C_3^kY_{s,r}^{3-k}:{Z_{s,r}}^k:, \\
   Y_{s,s}&=X(s,x) .\\
   \end{aligned}
   \right.
  \end{equation}
  Applying Theorem \ref{t6.5} with $Y_{t,r}$ replacing $Y_r$ we have
$$
\mathbb{E}\int_{t}^{t+1}\|Y_{t,r}\|_{H^1}^2dr\lesssim 1+\mathbb{E}\|Y_{t,t}\|_{H^{-1}}^2= 1+\mathbb{E}\|X(t,x)\|_{H^{-1}}^2.
$$
Combining this with (\ref{4.15.2}) we deduce that for $\alpha\in (0,1)$,
$$
\mathbb{E}\int_t^{t+1}\|X(s,x)\|_{\mathcal{C}^{-\alpha}}^2ds\leq \mathbb{E}\int_t^{t+1}\|Y_{t,s}\|_{H^1}^2ds+\mathbb{E}\int_t^{t+1}\|Z_{t,s}\|_{\mathcal{C}^{-\alpha}}^2ds\lesssim 1+\frac{1}{1\wedge t},
$$
where we used a similar argument as in the proof of \cite[Theorem 2.1]{Tsatsoulis:2016wm} in the last inequality.
Then we obtain that for $t\geq 1$
$$
\mathbb{E}\int_1^t\|X(s,x)\|_{\mathcal{C}^{-\alpha}}^2ds\lesssim t.
$$
Moreover, by (\ref{4.4}) we have
$$
\mathbb{E}\int_0^1\|Y_s\|_{H^{1}}^2ds\lesssim 1+\|x\|_{H^{-1}}^2.
$$
Thus for $t\geq 1$
$$
\int_0^t\mathbb{E}\|X(s,x)\|_{\mathcal{C}^{-\alpha}}^2ds\leq \int_0^1\mathbb{E}\|X(s,x)\|_{\mathcal{C}^{-\alpha}}^2ds+\int_1^t\mathbb{E}\|X(s,x)\|_{\mathcal{C}^{-\alpha}}^2ds\lesssim 1+\|x\|_{H^{-1}}^2+t.
$$
By Chebeshev's inequality, for any $K>0$
$$
\mathbb{P}(\|X(t,x)\|_{\mathcal{C}^{-\alpha}}>K)\leq\frac{1}{K^2}\mathbb{E}\|X(t,x)\|_{\mathcal{C}^{-\alpha}}.
$$
Thus there exists a constant $C>0$, such that
\begin{align*}
\int_0^t\mathbb{P}(\|X(s,x)\|_{\mathcal{C}^{-\alpha}}>K)&\leq \frac{C}{K^2}\int_0^t\mathbb{E}\|X(s,x)\|_{\mathcal{C}^{-\alpha}}^2ds\\
&\leq \frac{C}{K^2}(1+\|x\|_{H^{-1}}^2+t).
\end{align*}
Letting $R_t(x,\cdot)=\frac{1}{t}\int_0^tP_s(x,\cdot)ds$, for $K_{\varepsilon}^2=\frac{C}{\varepsilon}$ we get
$$
R_t(f\in \mathcal{C}^{-\alpha}\cap V_0^1:\|f\|_{\mathcal{C}^{-\alpha}}>K_{\varepsilon})\leq R_t(f\in V_0^1:\|f\|_{\mathcal{C}^{-\alpha}}>K_{\varepsilon})\leq (1+\frac{1+\|x\|_{H^{-1}}}{t})\varepsilon.
$$
By \cite[Proposition 4.6]{Triebel:2006vl} we know that $\{f\in \mathcal{C}^{-\alpha}\cap V_0^1:\|f\|_{\mathcal{C}^{-\alpha}}>K_{\varepsilon}\}$ is a compact subset of $V_0^{-1}$ since the embedding $\mathcal{C}^{-\alpha}\subset V^{-1}$ is compact. This implies the tightness of $\{R_t\}_{t\geq 0}$ in $V_0^{-1}$. By the Krylov-Bogoliubov existence theorem (see \cite[Corollary 3.1.2]{DaPrato:1996dr}), there exists a sequence $t_k\nearrow \infty$ and a measure $\nu_x$ such that $R_{t_k}\to \nu_x$ weakly in $V_0^{-1}$ and $\nu_x$ is an invariant measure for the semigroup $\{P_t\}_{t\geq 0}$.

$\hfill\Box$

To prove the exponential mixing property, we make use of the irreducibiltiy of $\underline{Z}$ and a uniform estimate, which is slightly different from that in the proof of \cite[Theorem 6.3]{Tsatsoulis:2016wm}.
\begin{theorem}\label{6.6.3}
There exists a constant $\lambda\in (0,1)$ and $T_0>0$ such that 
$$
\|P_t(x)-P_t(y)\|_{TV}\leq 1-\lambda,
$$
for every $x,y\in V_0^{-1}$, $t\geq T_0+1$.
\end{theorem}
\proof From (\ref{4.14.2}) we know that for any fixed $r>0$, there exist $T_0,M>0$ which are independent of $\omega,x$, such that for any initial value $x\in V_0^{-1}$, we have that $\{\omega: \normmm{\underline{Z}}_{T_0}\leq M \}\subset \{\|Y(T_0)\|_{V_0^{-1}}<\frac{r}{2}\}\cap\{\|Z(T_0)\|_{V_0^{-1}}<\frac{r}{2}\}$.

By Theorem \ref{t6.5.2} for every $a\in (0,1)$ there exists $r\equiv r(a)>0$ such that for every $x,y\in\bar{B}_r(0)$ and $t\geq 1$
\begin{equation}\label{6.11.2}
\|P_t(x)-P_t(y)\|_{TV}\leq 1-a,
\end{equation}
where $B_r(u):=\{x\in V_0^{-1}:\|x-u\|_{V_0^{-1}}<r\}$. Then by (\ref{4.14.2}) for any intial value $x\in V_0^{-1}$, there exists $b\equiv b(r)\in (0,1)$ such that
\begin{equation}\label{6.12.2}
\begin{aligned}
\mathbb{P}(\|X(T_0)\|_{V_0^{-1}}\leq r)&\geq\mathbb{P}\left(\{\|Y(T_0)\|_{V_0^{-1}}\leq \frac{r}{2}\}\cap\{\|Z(T_0)\|_{V_0^{-1}}\leq \frac{r}{2}\}\right)\\
&\geq \mathbb{P}(\normmm{\underline{Z}}_{T_0}\leq M)\\
&\geq b,
\end{aligned}
\end{equation}
where in the last step we used the irreducibility of the law of $\underline{Z}$. Here we omit the proof of the irreducibility of $\underline{Z}$, since it is the same as that of \cite[Theorem 6.3]{Tsatsoulis:2016wm}. Moreover, by (\ref{6.12.2}) for any $R>0$
\begin{equation}\label{6.13.2}
\inf_{x\in V_0^{-1}}P_{T_0}(x,\bar{B}_r(0))\geq b.
\end{equation}
Then combining (\ref{6.11.2})-(\ref{6.13.2}) and the Markov property by the same argument as in the proof of \cite[Theorem 6.5]{Tsatsoulis:2016wm}, we can complete the proof.
$\hfill\Box$
\vskip.10in

The following corollary gives the exponential convergence to a unique invariant measure.
\begin{corollary}
There exists a unique invariant measure $\bar{\nu}$ for the semigroup $\{P_t\}_{t\geq 0}$ such that
$$
\|P_t-\bar{\nu}\|_{TV}\leq(1-\lambda)^{\lfloor\frac{t}{T}\rfloor}\|\delta_x-\bar{\nu}\|_{TV}.,
$$
for every $x\in V_0^{-1}$, $t\geq T_0+1$. Moreover, $\bar{\nu}=\nu$.
\end{corollary}
\proof For the first result, see the proof of \cite[Corollary 6.6]{Tsatsoulis:2016wm}. By Corollary \ref{c5.7}, $\nu$ is an invariant measure of $X$.
 Hence $\bar{\nu}=\nu$.
$\hfill\Box$
\vskip.10in
\begin{remark}\label{r6.9}
In the following we give a simple and short proof for exponential convergence by the theory of Dirichlet forms. 

Similarly to \cite{DaPrato:2004fw}, by comparing the two Dirichlet forms for the Cahn-Hilliard equation and the dynamical $\Phi_2^4$ model, we obtain the spectral gap of equation (\ref{1.1}). Indeed, 
since $\nu$ is the restriction of $\Phi_2^4$ measure $\tilde{\nu}$ on $E$,
by the same arguments in \cite{Rockner:2015uh} and \cite{Tsatsoulis:2016wm} we know that $\tilde{\nu}$ is the invariant measure for the solution to the dynamical $\Phi_2^4$ model. We denote the Dirichlet form associated with the dynamical $\Phi_2^4$ model by $(\bar{\mathcal{E}},D(\bar{\mathcal{E}}))$, i.e. 
 $$
\bar{\mathcal{E}}(f,g)=\frac{1}{2}\int_{\tilde{E}}\langle Df,Dg\rangle_{L^2}d\nu, f,g\in D(\bar{\mathcal{E}}), 
 $$ 
 where $\tilde{E}=H^{-1-s_0}$,
  $D$ denotes the gradient operator in $L^2(\mathbb{T}^2)$ (see \cite{Rockner:2015uh}).
In \cite{Tsatsoulis:2016wm} the expoential convergence for the dynamical $\Phi_2^4$ model in total variation is proved. This implies the expoential convergence in $L^2(\tilde{E},\tilde{\nu)}$-norm. By \cite[Theory 1.1, Example 1.1.2]{Wang:2006uy} this is equivalent to the Poincar\'e inequality
$$
\int f^2d\tilde{\nu}-(\int fd\tilde{nu})^2\leq C\bar{\mathcal{E}}(f,f), f\in D(\bar{\mathcal{E}}).
$$
From the proof of Theorem \ref{t6.2} we know that for any $f\in D(\mathcal{E})\subset D(\bar{\mathcal{E}})$,
\begin{align*}
\mathcal{E}(f,f)=\frac{1}{2}\sum_{k\neq(0,0)}\int |\frac{\partial f}{\partial h_k}|^2d\nu
=\frac{1}{2}\sum_{k\neq(0,0)}\lambda_k \int |\frac{\partial f}{\partial e_k}|^2d\nu\geq\frac{1}{2}\sum_{k\neq(0,0)}\int |\frac{\partial f}{\partial e_k}|^2d\nu=\mathcal{\bar{E}}(f,f),
\end{align*}
where $h_k=\sqrt{\lambda_k}e_k$, $\{h_k\}_{k\in \mathbb{Z}^2\setminus \{(0,0)\}}$ is an orthonormal basis of $V_0^{-1}$. Then by \cite[Theory 1.1, Example 1.1.2]{Wang:2006uy} we have
$$
\|P_tf-\int fd\nu\|_{L^2(E,\nu)}\leq e^{-\frac{t}{C}}\|f-\int fd\nu\|_{L^2(E,\nu)}.
$$
  \end{remark}

   \appendix
  \renewcommand{\appendixname}{Appendix~\Alph{section}}
  
\section{Symmetric quasi regular Dirichlet forms and Markov Processes}\label{c2s2}

\vskip.10in
In this section we recall some general Dirichlet form results from \cite{Fukushima:1994kz,Ma:1992ec}.

 Let $E$ be a Hausdorff topological space, $m$ a $\sigma$-finite measure on $E$,
and let $\mathcal{B}$ the smallest $\sigma$-algebra of subsets of $E$ with respect to which all
continuous functions on $E$ are measurable.
Let $\mathcal{E}$ be a symmetric Dirichlet form acting in the real $L^2(m)$-space, i.e.
$\mathcal{E}$ is a positive, symmetric, bilinear, closed form with domain $D(\mathcal{E})$ dense
in $L^2(m)$, and such that $\mathcal{E}(\Phi(u),\Phi(u))\leq \mathcal{E}(u, u)$, for any $u\in D(\mathcal{E})$, where
$\Phi(t) = (0\vee t)\wedge1, t\in\mathbb{R}$. The latter condition is known to be equivalent
with the condition that the associated $C_0$-contraction semigroup $T_t, t\geq 0$,
is submarkovian (i.e. $0 \leq u\leq 1$ m-a.e. implies $0\leq T_tu \leq 1$ m-a.e., for
all $u\in L^2(m)$); association means that $\lim_{t\downarrow0}\frac{1}{t}\langle u- T_tu, v\rangle_{L^2(m)} = \mathcal{E}(u, v), \forall u, v\in D(\mathcal{E})$.
\vskip.10in

 \begin{definition}(\cite[Chapter \Rmnum{3}, Definition 2.1, Definition 3.2]{Ma:1992ec})\label{dqe}
  
   (i) An increasing sequence $(F_k)_{k\in\mathbb{N}}$ of closed subsets of $E$ is called an $\mathcal{E}$-\emph{nest} if $\cup_k\{f\in D(\mathcal{E}):\operatorname{supp}f\subset F_k\}$ is dense in $D(\mathcal{E})$ w.r.t. $\mathcal{E}_1^{\frac{1}{2}}$. Here $\mathcal{E}_1^{\frac{1}{2}}$ is a norm on $E$
 defined by $\mathcal{E}_1(u):=\mathcal{E}(u,u)+\|u\|_{L^2(E,\nu)}^2$, $\forall u\in E$.   
 
   (ii) A subset $N\subset E$ is called $\mathcal{E}$-\emph{exceptional} if $N\subset \cap_{k}F_k^c$ for some $\mathcal{E}$-nest $(F_k)_{k\in\mathbb{N}}$. We say that a property of points in $E$ holds $\mathcal{E}$-\emph{quasi everywhere}(abbreviated $\mathcal{E}-q.e.$) if the property holds outside some $\mathcal{E}$-exceptional set.
   
   (iii) An $\mathcal{E}-q.e.$ defined function $f$ on $E$ is called $\mathcal{E}$-\emph{quasi-continuous} if there exists an $\mathcal{E}$-nest $(F_k)_{k\in\mathbb{N}}$ such that $f\big|_{F_k}$ is continuous for every $k\in \mathbb{N}$.
   
   \end{definition}

\begin{definition}\label{qrdf}  (cf. \cite[Chap. \Rmnum{4}, Defi. 3.1]{Ma:1992ec})
A symmetric Dirichlet form is called \emph{quasi-regular} if the following holds:

(\rmnum{1}) There exists an $\mathcal{E}$-{nest} of $E$ such that $F_k$ is compact in $E$ for any $k$.

(\rmnum{2}) There exists an $\mathcal{E}_1^{1/2}$-dense subset of $D(\mathcal{E})$ whose elements have $\mathcal{E}$-quasi continuous
$m$-versions. A real function $u$ on $E$ is called quasi continuous
when there exists an $\mathcal{E}$-nest $(F_k)$ s.t. $u$ restricted to $F_k$ is continuous.

(\rmnum{3}) There exists $u_n\in D(\mathcal{E}), n\in\mathbb{N}$, with $\mathcal{E}$-quasi continuous $m$-versions $˜\tilde{u}_n$
and there exists an $\mathcal{E}$-exceptional subset $N$ of $E$ s.t. $\{˜\tilde{u}_n\}_{n\in\mathbb{N}}$ separates
the points of $E\setminus N$. An $\mathcal{E}$-exceptional subset of $E$ is a subset
$N\subset \cap_k(E\setminus F_k)$ for some $\mathcal{E}$-nest $(F_k)$.
\end{definition}
\vskip.10in

To recall the main results in \cite{Ma:1992ec} we recall the definitions of a Markov process and a right process. Here we consider only Markov processes with life time $\infty$.
\vskip.10in

\begin{definition}\label{mp} (cf. \cite[Chap. \Rmnum{4} Defi. 1.5]{Ma:1992ec} A collection $\mathbf{M}:=(\Omega,\mathcal{M},(X_t)_{t\geq0},(\mathbb{P}^z)_{z\in E})$ is called a Markov process (with state space $E$) if
it has the following properties.

(\rmnum{1}) There exists a filtration $(\mathcal{M}_t)$ on $(\Omega,\mathcal{M})$ such that $(X_t)_{t\geq0}$ is an $(\mathcal{M}_t)_{t\geq0}$ adapted stochastic process
 with state space $E$.

(\rmnum{2}) For each $t\geq0$ there exists a shift operator $\theta_t:\Omega\rightarrow\Omega$ such that $X_s\circ\theta_t=X_{s+t}$ for all $s,t\geq0$

(\rmnum{3}) $\mathbb{P}^z, z\in E,$ are probability measures on $(\Omega,\mathcal{M})$ such that $z\mapsto \mathbb{P}^z(A)$ is $\mathcal{B}(E)^*$-measurable for each $A\in\mathcal{M}$ resp. $\mathcal{B}(E)$-measurable if $A\in \sigma\{X_s|s\in[0,\infty)\}$, where $\mathcal{B}(E)^*:=\cap_{\mathbb{P}\in\mathcal{P}(E)}\mathcal{B}^\mathbb{P}(E)$ for $\mathcal{P}(E)$ denoting the family of all probability measures on $(E,\mathcal{B}(E))$ and $\mathcal{B}^\mathbb{P}(E)$ denotes the completion of the $\sigma$-algebra $\mathcal{B}(E)$ w.r.t. a probability $\mathbb{P}$.

(\rmnum{4}) (Markov property) For all $A\in\mathcal{B}(E)$ and any $t,s\geq0$
$$\mathbb{P}^z[X_{s+t}\in A|\mathcal{M}_s]=\mathbb{P}^{X_s}[X_t\in A] \quad \mathbb{P}^z-a.s., z\in E.$$
\end{definition}
\vskip.10in

\begin{definition} (cf. \cite[Chap. \Rmnum{4} Defi. 1.8]{Ma:1992ec}) Let $\mathbf{M}:=(\Omega,\mathcal{M},(X_t)_{t\geq0},(\mathbb{P}^z)_{z\in E})$ be a Markov process with state space $E$ and corresponding filtration $(\mathcal{M}_t)$. $\mathbf{M}$ is called a \emph{right process} if it has the following additional properties.

(\rmnum{1}) (Normal property) $\mathbb{P}^z(X_0=z)=1$ for all $z\in E$.

(\rmnum{2}) (Right continuity) For each $\omega\in\Omega$, $t\mapsto X_t(\omega)$ is right continuous on $[0,\infty)$.

(\rmnum{3}) (Strong Markov property) $(\mathcal{M}_t)$ is right continuous and for every $(\mathcal{M}_t)$-stopping time $\sigma$ and every $\nu\in \mathcal{P}(E)$
$$\mathbb{P}^\nu[X_{\sigma+t}\in A|\mathcal{M}_\sigma]=\mathbb{P}^{X_\sigma}[X_t\in A]\quad P^\nu-a.s.$$
for all $A\in\mathcal{B}(E)$, $t\geq0$.
\end{definition}
\vskip.10in

\begin{theorem}\label{tqrdf} (\cite[Chap. \Rmnum{4} Thm 6.7]{Ma:1992ec}) Let $E$ be a metrizable Lusin space. Then a Dirichlet form $(\mathcal{E},D(\mathcal{E}))$ on $L^2(E,m)$ is quasi-regular if and only if there exists a right process $\mathbf{M}$ associated with $(\mathcal{E},D(\mathcal{E}))$, i.e. the semigroup of $\mathbf{M}$ is an $m$-version of the semigroup associated with $(\mathcal{E},D(\mathcal{E}))$. In this case $\mathbf{M}$ is always properly associated with  $(\mathcal{E},D(\mathcal{E}))$.
\vskip.10in
\end{theorem}

\begin{remark}
 The results in \cite[Chap. \Rmnum{4}]{Ma:1992ec} are more general and can be applied for general Hausdorff topological spaces and more general Markov processes.
 Lusin spaces are enough for our use in this thesis.
 \end{remark}
 
Let us recall the definition of additive functional in \cite[Chapter 5]{Fukushima:1994kz}
 
 \begin{definition}\label{AF}
A real valued function $A_t(\omega)$, $t\geq 0$, $\omega\in \Omega$ is called an additive functional (AF in abbreviation) of right process $\mathbf{M}$ in Defition \ref{mp} if 

(A.1) $A_t$ is $\mathcal{F}_t$-measurable, where $\{\mathcal{F}_t\}$ is the natural filtration of $\mathbf{M}$.

(A.2) There exist a set \(\Lambda \in \mathcal{F}_\infty\) and an exceptional set \(N \subset E\) such that \(\mathbb{P}^{x}(\Lambda)=1\), \(\forall x \in E \backslash N, \theta_{t} \Lambda \subset \Lambda, \forall t>0\), and moreover, for each \(\omega \in \Lambda\),
 \(A_.(\omega)\) is right continuous and has the left limit on \([0, \infty), A_{0}(\omega)=0,\left|A_{t}(\omega)\right|<\infty, \forall t<\infty$ and \(A_{t+s}(\omega)=\)
\(A_{s}(\omega)+A_{t}\left(\theta_{s} \omega\right), \forall t, s \geq 0\).
 \end{definition}

  \section*{Acknowledgements}
  We are very grateful to Professor Zhi-Ming Ma and Professor Xiangchan Zhu for numerous discussions.

\bibliographystyle{alpha}%

\end{document}